\documentclass[11pt]{amsart}

\usepackage{lmodern}
\usepackage{textcomp}
\usepackage{epic}

\usepackage{enumitem}

\usepackage{amsmath,amsthm,amsfonts,amssymb,mathtools,mathrsfs,bm}

\usepackage{aliascnt}

\usepackage{accents}

\usepackage{nicefrac}

\usepackage{tikz-cd}
\usetikzlibrary{babel}
\usetikzlibrary{decorations.pathmorphing}

\usepackage[utf8]{inputenc}

\usepackage[T1]{fontenc}

 \usepackage{lmodern}

\usepackage{float}

\usepackage{xcolor}
\definecolor{colordelink}{rgb}{0,0,0.50}
\definecolor{colordecite}{rgb}{0,0.5,0}
\definecolor{colordeurl}{rgb}{0,0.41,0.5}

\usepackage[pagebackref=true]{hyperref}
\hypersetup{
    colorlinks=true,
    linkcolor=colordelink,
    citecolor=colordecite,
	 urlcolor=colordeurl,
	 final
}

\usepackage[capitalize,nameinlink,noabbrev]{cleveref}

\usepackage{todonotes}

\usepackage{multirow}
\usepackage{cellspace}
\def\dim{\operatorname{dim}}

\def\codim{\operatorname{codim}}

\newcommand{\CC}{\mathbb{C}}

\newcommand{\eqA}{\mathscr{A}}

\newcommand{\Ascr}{\mathscr{A}}

\newcommand{\Kscr}{\mathscr{K}}

\newcommand{\Rscr}{\mathscr{R}}

\newcommand{\Nbb}{\mathbb{N}}

\newcommand{\Mcal}{\mathcal{M}}

\newcommand{\Ocal}{\mathcal{O}}

\newcommand{\ubm}{\bm{u}}
\newcommand{\vbm}{\bm{v}}
\newcommand{\wbm}{\bm{w}}
\newcommand{\xbm}{\bm{x}}
\newcommand{\ybm}{\bm{y}}

\makeatletter
\newcommand{\tpitchfork}{%
  \vbox{
    \baselineskip\z@skip
    \lineskip-.52ex
    \lineskiplimit\maxdimen
    \m@th
    \ialign{##\crcr\hidewidth\smash{$-$}\hidewidth\crcr$\pitchfork$\crcr}
  }%
}
\makeatother

\theoremstyle{plain}
\newtheorem{theorem}{Theorem}[section]

\newaliascnt{lemma}{theorem}
\newtheorem{lemma}[lemma]{Lemma}
\aliascntresetthe{lemma}

\newaliascnt{corollary}{theorem}
\newtheorem{corollary}[corollary]{Corollary}
\aliascntresetthe{corollary}

\newaliascnt{proposition}{theorem}
\newtheorem{proposition}[proposition]{Proposition}
\aliascntresetthe{proposition}

\newaliascnt{question}{theorem}

\aliascntresetthe{question}

\newaliascnt{thm}{theorem}

\aliascntresetthe{thm}

\newaliascnt{lem}{theorem}

\aliascntresetthe{lem}

\newaliascnt{cor}{theorem}

\aliascntresetthe{cor}

\newaliascnt{prop}{theorem}

\aliascntresetthe{prop}
\Crefname{lemma}{Lemma}{Lemmas}
\crefname{lemma}{Lemma}{Lemmas}
\Crefname{proposition}{Proposition}{Propositions}
\crefname{proposition}{Proposition}{Propositions}
\Crefname{corollary}{Corollary}{Corollaries}
\crefname{corollary}{Corollary}{Corollaries}

\newaliascnt{diagram}{equation}
\aliascntresetthe{diagram}
\crefname{diagram}{Diagram}{Diagrams}
\Crefname{diagram}{Diagram}{Diagrams}

\Crefname{prop}{Proposition}{Propositions}
\crefname{prop}{Proposition}{Propositions}
\Crefname{lem}{Lemma}{Lemmas}
\crefname{lem}{Lemma}{Lemmas}
\Crefname{cor}{Corollary}{Corollaries}
\crefname{cor}{Corollary}{Corollaries}

\Crefname{definition}{Definition}{Definitions}
\crefname{definition}{Definition}{Definitions}

\newtheorem*{theorem**}{Theorem\theoremnum}
\newenvironment{theorem*}[1][]{%
  \edef\theoremnum{\if\relax\detokenize{#1}\relax\else~#1\fi}
  \begin{theorem**}
}{%
  \end{theorem**}
} 
\newtheorem*{conjecture**}{Conjecture\theoremnum}
\newenvironment{conjecture*}[1][]{%
  \edef\theoremnum{\if\relax\detokenize{#1}\relax\else~#1\fi}
  \begin{conjecture**}
}{%
  \end{conjecture**}
} 
\newtheorem*{corollary**}{Corollary\theoremnum}
\newenvironment{corollary*}[1][]{%
  \edef\theoremnum{\if\relax\detokenize{#1}\relax\else~#1\fi}
  \begin{corollary**}
}{%
  \end{corollary**}
}

\theoremstyle{definition}
\newaliascnt{definition}{theorem}
\newtheorem{definition}[definition]{Definition}
\aliascntresetthe{definition}

\newaliascnt{proposition/definition}{theorem}
\newtheorem{proposition/definition}[proposition/definition]{Proposition/Definition}
\aliascntresetthe{proposition/definition}

\newaliascnt{conjecture}{theorem}

\aliascntresetthe{conjecture}

\newaliascnt{example}{theorem}
\newtheorem{example}[example]{Example}
\aliascntresetthe{example}

\newaliascnt{oproblem}{theorem}

\aliascntresetthe{oproblem}

\newaliascnt{remark}{theorem}
\newtheorem{remark}[remark]{Remark}
\aliascntresetthe{remark}

\newaliascnt{notation}{theorem}
\newtheorem{notation}[notation]{Notation}
\aliascntresetthe{notation}
\Crefname{notation}{Notation}{Notations}

\theoremstyle{remark}

\begin{document}
\author{R. Giménez Conejero}
\address{Escuela Técnica Superior de Ingeniería y Diseño Industrial (ETSIDI), Universidad Politécnica de Madrid, Ronda de Valencia 3, 28012 Madrid, Spain}
\email{Roberto.Gimenez@uv.es}

\author{M. E. Rodrigues Hernandes}
\address{Departamento de Matemática, Universidade Estadual de Maringá, Av.~Colombo 5790, 87020-900 Maringá, PR, Brazil}
\email{merhernandes@uem.br}

\title[Monomial singularities of maps in low codimension]{Monomial singularities of maps in low codimension and their deformations}


\begin{abstract}
We study complex analytic map germs whose coordinate functions are monomials up to a change of coordinates. Motivated by recent normal form results of Rodrigues Hernandes and Ruas for $\Ascr$-finite monomial maps in the range $p\ge 2n$, we address the remaining gap $p<2n$, solving an open problem posed in their work in the range $n<p<2n$. Our main results show that a monomial map germ $f:(\CC^n,0)\to(\CC^p,0)$ with $p < 2n$ and $(n,p) \neq (1,1)$ is $\Ascr$-finite if and only if it is $\Ascr$-equivalent to an $\Ascr$-stable germ. For $n = p = 1$ we also obtain the family of simple singularities of type $A_k$.

Finally, we construct new families of $\Ascr$-finite (generally non-monomial) map germs starting from monomial ones.

\end{abstract}

\maketitle

\section{Introduction}

In this work we focus on maps whose coordinate functions are monomials, modulo $\Ascr$-equivalence (i.e., change of coordinates). Since these maps are quasi-homogeneous, considering their global behavior or their germs at the origin is essentially equivalent. We follow the latter approach, as commonly done in singularity theory.

Monomial map germs play a fundamental role in Singularity Theory and Commutative Algebra. The interest in studying this class of germs is not merely justified by their operational simplicity and combinatorial appeal; for instance, they provide concrete analytic models and parametrize toric varieties. Recall that an affine toric variety can be defined as an algebraic variety whose defining ideal is prime and generated by binomials. Moreover, the exponents of the monomial parametrization of a toric variety directly correspond to the generators of the combinatorics of this variety (i.e., its semigroup). Consequently, studying the properties of monomial map germs is equivalent to decoding the structure of toric singularities through the combinatorial properties of their exponent vectors. See, for instance, \cite{Dalbelo2026,Araujo2026}.

From the perspective of Singularity Theory, key properties such as $\Ascr$-finite determinacy (see \cref{def:a-finite}), the computation of analytic invariants, and the study of deformations are more manageable for monomial classes than for general map germs. 

The property of $\Ascr$-finite determinacy is essential in the local theory of smooth maps; for instance, it was shown in \cite{Bobadilla2019} that the Mond conjecture, one of the central open problems in singularities of maps, holds if it is verified for a suitable family of arbitrary multiplicity. Since monomial map germs are generally not $\mathscr{A}$-finite, this result cannot be applied directly to such classes.

Such considerations motivated the study of $\Ascr$-finite monomial map germs by M.E. Rodrigues Hernandes and M.A.S. Ruas. In \cite{RodriguesHernandes2019}, they provided normal forms for monomial surfaces parametrized by $\Ascr$-finite monomial map germs $f:(\CC^2,0)\to(\CC^4,0)$. In \cite{RodriguesHernandes2023}, they extended their results to establish a general form for $\Ascr$-finite monomial germs $f:(\CC^n,0)\to(\CC^p,0)$ with $p\geq 2n$. To construct such families, they introduced the notion of a \textit{join operation}, which relies on analytic invariants of these germs, such as the delta-invariant and the semigroup associated with a certain parametrized curve. The delta-invariant generalizes the classical invariant for curves and is easier to compute than the $\Ascr_e$-codimension. In particular, they obtained some normal forms, as exemplified in the following result when $p=2n$: 
\begin{theorem}[see {\cite[Corollary 4.3]{RodriguesHernandes2023}}]
    A monomial map germ $g:(\CC^n,0)\to(\CC^{2n},0)$ is $\Ascr$-finite if and only if it is $\Ascr$-equivalent to an immersion or to
    \[f(\xbm,y)=(\xbm,y^{m_1},y^{m_2},x_1^{\lambda_1}y,\dots,x_{n-1}^{\lambda_{n-1}}y),\]
    where $\lambda_i\geq 1$ for all $i=1,\dots,n-1$, and $\gcd(m_1,m_2)=1$.
\end{theorem}

 Therefore, with respect to $\Ascr$-finite monomial maps, the work \cite{RodriguesHernandes2023} leaves unanswered the dimensions $n < p<2n$ and $n \geq p$. Indeed, in \cite[Remark 4.8]{RodriguesHernandes2023}, the authors observed that all $\Ascr$-finite monomial map germs $f:(\CC^n,0)\to(\CC^{p},0)$, with $n<p<2n$, obtained  by their approach are immersions or generalized cross-caps, implicitly conjecturing that these are the only possible classes. The proof of this conjecture is one of our main results.

\begin{theorem*}[A]\hypertarget{thm:monomialAfinite}{}
Let $f:(\CC^n,0)\to(\CC^{p},0)$ be a monomial map germ, for $n < p < 2n$. Then $f$ is $\Ascr$-finite if and only if $f$ is $\Ascr$-equivalent to one of the following map germs:
\begin{enumerate}
    \item $(\xbm,y)\mapsto(x_1,\ldots, x_{n-1},y,0, \ldots, 0)$,
    \item $(\xbm,y)\mapsto(x_1,\ldots, x_{n-1}, y^2,  x_{1}y,\ldots,x_{p-n}y)$,
\end{enumerate} 
where $\xbm :=(x_1,\ldots, x_{n-1})$, and all are $\Ascr$-stable.
\end{theorem*}

To complete the classification of $\Ascr$-finite monomial map germs for all pairs of dimensions $(n,p)$, we show that for $p \leq n$ with $(n,p) \neq (1,1)$, all such germs are $\Ascr$-stable. Finally, for $n = p = 1$ we also obtain the family of simple singularities of type $A_{a-1}$ ($a \geq 2$).

\begin{theorem*}[B]\hypertarget{thm:ngeqp}{}
Let $f:(\CC^n,0)\to(\CC^{p},0)$ be a monomial map germ, for $p\leq n$. Then $f$ is $\Ascr$-finite if and only if $f$ is $\Ascr$-equivalent to one of the following map germs:
\begin{enumerate}
    \item $(x_1,\dots,x_n)\mapsto(x_1,\dots,x_p)$,
    \item $(x_1,\dots,x_n)\mapsto(x_1,\dots,x_{n-1},x_n^2)$, if $n=p$,
    \item $(x_1,\dots,x_n)\mapsto(x_1,\dots,x_{p-1},x_px_{p+1})$, if $n=p+1$,
    \item $x\mapsto x^a$, with $a\geq 1$, if $n=p=1$.
\end{enumerate}
The germs in items (1)-(3) are $\Ascr$-stable, and the germ in item (4) is $\Ascr$-stable for $a=1,2$ and for $a> 2$ has $\Ascr_e$-codimension $a-2$.
\end{theorem*}

In the classification work \cite{RodriguesHernandes2023}, as well as in \hyperlink{thm:monomialAfinite}{Theorem A} and \hyperlink{thm:ngeqp}{Theorem B}, some normal forms are given for the $\Ascr$-classes of monomial map germs. However, there might be redundancies because two normal forms may still be $\Ascr$-equivalent. In the recent work \cite{In_prep}, R. Giménez Conejero and D. Gonzalez-Sanchez show that two monomial maps (without any assumption on $\Ascr$-finiteness or on the dimensions) are $\Ascr$-equivalent if and only if their exponent semigroups coincide up to permutation of variables. In particular, it shows that two different normal forms are not $\Ascr$-equivalent, completely closing the classification of $\Ascr$-finite monomial map germs. 
\newline

We motivate the proof of \hyperlink{thm:monomialAfinite}{Theorem A} with examples in \cref{sec:Examples}, and we give the proof in \cref{sec:corankone,sec:coranktwo}. We present two proofs of \hyperlink{thm:monomialAfinite}{Theorem A} for corank one map germs, relying on distinct techniques introduced in \cref{sec:Preliminaries}. The first employs a geometric approach via the Marar-Mond criterion, where the multiple point scheme $D^2(f)$ translates stability into an explicit ICIS verification at the double point locus. The second is purely algebraic and uses the finiteness of the $\mathcal{A}_e$-codimension, providing a framework that extends to higher corank germs.

The proof of \hyperlink{thm:ngeqp}{Theorem B} is given in \cref{sec:ngep}, using similar algebraic ideas.

These techniques are also applied in \cref{sec:Non-monomial}, where we construct new $\mathscr{A}$-finite germs starting from monomial maps.
\newline

\textbf{Acknowledgements:} The authors want to thank the organizers of the \textit{Brazilian Singularity Theory Webinar}: Débora Lopes, 
Hellen Santana,
Igor Chagas, 
Inácio Rabelo,
Jorge L. Deolindo-Silva, and Raphael de Oliveira. This work originated after a talk by the second-named author.

The authors also thank Ignacio Breva Ribes for his comments on augmentations of singularities, and Juan J. Nuño-Ballesteros for providing a proof of Theorem A in the case of surfaces in $\mathbb{C}^3$.

\section{Preliminaries}\label{sec:Preliminaries}

\subsection{Map germs singularities}
In this subsection, we introduce some basic concepts of singularities of maps. For a general reference, we suggest \cite{Mond2020}.

\begin{notation}
The set of smooth map germs $f:(\CC^n,0)
\to \CC^p$ is denoted by $\mathcal{O}(n,p)$. The local ring $\mathcal{O}_n\coloneqq \mathcal{O}(n,1)$ has maximal ideal $\mathcal{M}_n=\{f \in \mathcal{O}_n; f(0)=0\}$. 

Observe that $\mathcal{O}(n,p)$ is a free $\mathcal{O}_n$-module of rank $p$, and the submodule of germs $f :(\CC^n,0)\to(\CC^p,0)$
is denoted by $\mathcal{M}_n \mathcal{O}(n,p)$.

We will also need to consider the
set of germs of vector fields in $(\CC^k,0)$, denoted by $\theta_k$, and the set $\theta(f)$ of germs of vector fields along $f$ defined as $\theta(f)=\{\sigma: (\CC^n,0)\to T\CC^p; \ \pi \circ \sigma =f\}$, where $\pi: T\CC^p \to (\CC^p,0)$ is the natural projection.
\end{notation}

 The classical group of change of coordinates in
source and target, denoted by $\Ascr$, is the set of germs of analytic
diffeomorphisms $(\rho, \sigma)$ where $\rho:(\CC^n, 0) \to
(\CC^{n},0)$ and $\sigma:(\CC^p, 0) \to (\CC^{p},0)$. An element
$(\rho, \sigma) \in \Ascr$ acts on $f \in \mathcal{M}_n\mathcal{O}(n,p)$ as $(\rho,\sigma)\cdot f =\sigma \circ f \circ \rho^{-1}$. Two germs
$f$ and $g$ are called \textit{$\Ascr$-equivalent} if they are in the same $\Ascr$-orbit.

\begin{definition}
    The {\it extended tangent space} to the $\Ascr$-orbit of $f\in \Mcal_n\Ocal(n,p)$ is
\[T\Ascr_e(f)=tf(\theta_n) + wf(\theta_p),\]
 where
\[
\begin{aligned}
tf &: \theta_n \to \theta(f), & tf(\epsilon) &= df(\epsilon), \\
wf &: \theta_p \to \theta(f), & wf(\eta) &= \eta \circ f,
\end{aligned}
\] where $df$ is the differential of $f$.
Finally, the $\Ascr_e$-codimension of $f\in\Mcal_n\Ocal(n,p)$ is the dimension as $\CC$-vector space of the $\Ascr_e$-normal space,
\[N\Ascr_e(f)\coloneqq\frac{\theta(f)}{T\Ascr_e(f)},\]
and it is denoted by
\[\codim_{\Ascr_e}(f)\coloneqq \dim_{\CC}N\Ascr_e(f).\]
\end{definition}

A map germ $f$ is called {\it $\Ascr$-stable} (or simply stable) if $\codim_{\Ascr_e}(f)=0$.

In coordinates, for $f(\xbm)=\big(f_1(\xbm),\ldots, f_p(\xbm)\big)$ where $\xbm=(x_1,\ldots, x_n)$, we write
\begin{equation}
   \label{Tang-Space}
   T\Ascr_e(f)=\left\{\left(\begin{array}{ccc}
\frac{\partial f_1}{\partial x_1} & \cdots & \frac{\partial f_1}{\partial x_n}\\
\vdots & \vdots & \vdots \\
\frac{\partial f_p}{\partial x_1} & \ldots & \frac{\partial f_p}{\partial x_n}\\
\end{array}\right) \left(\begin{array}{c}
\epsilon_1\\
\vdots\\
\epsilon_n\\
\end{array}\right)+\left(\begin{array}{c}
\eta_1(f)\\
\vdots\\
\eta_p(f)
\end{array}\right); \begin{array}{c}
\epsilon \in \theta_n, \ \eta \in \theta_p,\\
\epsilon=(\epsilon_1,\ldots,\epsilon_n) \\
\eta=(\eta_1,\ldots,\eta_p)
\end{array}\right\}.
\end{equation} To simplify notation, we can represent the elements in $T\Ascr_e(f)$ as vectors in $\mathbb{C}^p$. We also denote by $e_i\in\theta(f)$ the canonical vector that has zeros in all components except at the $i$-th component, which is one. 

Note that 
$$wf(\theta_p)=f^*\left(\bigoplus_{i=1}^p \Ocal_{p}\right)\simeq \bigoplus_{i=1}^p f^*(\Ocal_p),$$ where $f^*(\Ocal_p)=\{g \circ f; \ g \in \Ocal_p \}$.

\begin{definition}\label{def:a-finite}
    A smooth map germ $f\in\Mcal_n\Ocal(n,p)$ is called \textit{$\Ascr$-finitely determined} (or \textit{${\Ascr}$-finite}) if there exists $k < \infty$ such that for all map germs $g\in\Mcal_n\Ocal(n,p)$ with $j^kf=j^kg$, then $f$ is ${\Ascr}$-equivalent to $g$, where $j^kf$ denotes the $k$-jet of $f$ at the origin.
\end{definition}

By Mather's Infinitesimal criterion of finite determinacy, \cite{Mather1968}, $f$ is ${\Ascr}$-finite if and only if $\codim_{\Ascr_e}(f)$ is finite. 

\subsection{Multiple point spaces}

For a map germ $f:(\CC^n,0)\rightarrow (\CC^p,0)$ with $n < p$, the multiple point spaces $D^k(f)$ encode the points in the image of $f$ that have $k$-preimages (or, rather, their closure). A good introduction to these topics is provided in \cite[Chapter 2]{Robertothesis} and \cite{Houston1999}.

Let $f:(\CC^n,0)\rightarrow (\CC^p,0)$ be a corank one map germ with $n< p$, which can be assumed, up to a change of coordinates, to be given by
\begin{equation}\label{Corank-1-germ}
f(\xbm,y)=\big(\xbm,f_1(\xbm,y)\ldots,f_{p-n+1}(\xbm,y)\big), \end{equation} where $\xbm=(x_1,\ldots, x_{n-1})$.

The $k$-th multiple point space $D^k(f)$ is defined as an analytic subvariety of the ambient space $(\mathbb{C}^{n+k-1}, 0)$ with coordinates $(\xbm, y_1, \dots, y_k)$. The space $D^1(f)$ is simply the source of the germ, $(\CC^n,0)$. For $k>1$ and $p-k(p-n)\geq 0$, a set of defining equations for $D^k(f)$ is given by the  \textit{divided differences} $f_j^{(i+1)}(\xbm, y_1, \ldots, y_{i+1}) = 0$, for $1 \le j \le p-n+1$ and $1 \le i \le k-1$, defined inductively by 
\[\begin{aligned}
f_j^{(1)}(\xbm,y_1)\coloneqq &f_j(\xbm,y_1),  \\
f_j^{(i+1)}(\xbm,y_1,\ldots,y_{i+1}) \coloneqq&\frac{f_j^{(i)}(\xbm,y_1,\ldots,y_{i-1},y_{i})-f_j^{(i)}(\xbm,y_1,\ldots,y_{i-1},y_{i+1})}{y_{i}-y_{i+1}}, \   \mbox{for} \ i\geq 1.
\end{aligned}\] We illustrate these concepts in \cref{ex:examplestableDk}.

 The space $D^k(f)$ has ambient space $(\CC^{n+k-1},0)$, and a set of defining equations is given by the  \textit{divided differences}:
\[\begin{aligned}f_j^{(i+1)}(\xbm,y_1,\ldots,y_{i+1}) \coloneqq&\frac{f_j^{(i)}(\xbm,y_1,\ldots,y_{i-1},y_{i})-f_j^{(i)}(\xbm,y_1,\ldots,y_{i-1},y_{i+1})}{y_{i}-y_{i+1}},\\ f_j^{(1)}(\xbm,y_1)\coloneqq&f_j(\xbm,y_1),
\end{aligned}\] with $1 \leq j \leq p-n+1$.
The space $D^1(f)$ is simply the source of the germ, $(\CC^n,0)$. The equations of the space $D^k(f)$, when $k>1$ and $p-k(p-n)\geq 0$ are $f_j^{(i+1)}$  $j=1,\dots,p-n+1$ and $i=1,\dots,k-1$. This is illustrated in \cref{ex:examplestableDk}.

Let us expand on what the divided differences are. 

\begin{definition}\label{def:Homog-Polyn-d}
 For each $d \ge 0$, let $h_d$ denote the \textit{complete homogeneous symmetric polynomial} of degree $d$ in $n$ variables $x_1,\dots,x_n$, i.e.,
\[h_d(x_1,\dots,x_n)=\sum_{1\leq i_{1}\leq i_{2}\leq \cdots \leq i_{d}\leq n}x_{i_{1}}x_{i_{2}}\cdots x_{i_{d}}.\]
If $d<0$ we simply set $h_d=0$ and $h_0(x_1, \dots, x_n)=1$.
The \textit{generating function} of these polynomials is
\begin{equation}\label{generatingfunction}
H_n(t)(x_1,\dots,x_n) \coloneqq  \prod_{i=1}^n \frac{1}{1 - x_i t}= \sum_{d \ge 0} h_d(\xbm) t^d.
\end{equation}
\end{definition}

These complete symmetric homogeneous polynomials are the building bricks of the divided differences.
Indeed, if we consider a function germ $u:(\CC^{k}\times\CC,0)\to(\CC,0)$ such that
\[u(\xbm,y)=\sum_j u_j(\xbm)y^j,\]
then it is easy to see that its divided difference is
\[u^{(2)}(\xbm,y_1,y_2):=\frac{u(\xbm,y_1)-u(\xbm,y_2)}{y_1-y_2}=\sum_j u_j(\xbm) h_{j-1}(y_1,y_2).\]
More generally, for the $k$-th iterated divided difference one has
\begin{equation}\label{eq:bricks}
u^{(k)}(\xbm,y_1,\dots,y_k)=\sum_ju_j(\xbm)h_{j-k+1}(y_1,\dots,y_k).
\end{equation}

\begin{remark}\label{rem:factorization}
Notice that, when $u$ is a monomial, the divided difference in \cref{eq:bricks} defines an irreducible variety (or empty) exactly when:
\begin{itemize}
\item $u(\xbm,y)=x_i^a y^{k-1}$, for some $i$ and any $k$;
\item $u(\xbm,y)= y^\alpha$ and $k>2$.
\item $u(\xbm,y)= y^\alpha$, $\alpha=1, 2$ and $k=2$.
\end{itemize}
\end{remark}


We conclude this section with a well-known evaluation of complete homogeneous symmetric polynomials, slightly adapted here to incorporate a scaling parameter $\lambda \in \mathbb{C}$. This property plays a key role in identifying points in the multiple point spaces. For the reader's convenience, we include its short proof below (see, for instance, \cite[Chapter I, Section 2]{MacDonald1998} and \cite[Corollary 1.5]{DeCarli2026}). 

\begin{proposition}
\label{thm:sympolZEROS}
Let $n, d \in \mathbb{N}$, $\lambda \in \mathbb{C}$ and $\zeta$ be a primitive $n$-th root of unity. Then \[
h_d(\lambda,\lambda\zeta,\dots,\lambda\zeta^{n-1})=
  \begin{cases}
    \lambda^{d}, & \text{if } n \mid d \\
    \ 0, & \text{if } n\nmid d .
  \end{cases}\qedhere
\] 
\end{proposition}
\begin{proof}
First, we show that
\begin{equation}\label{eq:productlambdas}
\prod_{j=0}^{n-1} (1 - \zeta^j t) = 1 -  t^n.
\end{equation}
Indeed, since $$\displaystyle\prod_{j=0}^{n-1} ( x - \zeta^j )=x^n-1,$$
evaluating at $x=t^{-1}$ yields
\[\prod_{j=0}^{n-1} \left(\frac{1}{t} - \zeta^j\right)= \left(\frac{1}{t}\right)^n-1.\]
\cref{eq:productlambdas} is obtained by multiplying the previous equality by $t^n$.

Now, by \cref{generatingfunction,eq:productlambdas}, we obtain
\[
H_n(t)(1,\zeta,\dots,\zeta^{n-1})=\prod_{j=1}^{n} \frac{1}{1 - \zeta^{j-1} t}
= \frac{1}{1 - t^n}
= \sum_{k \ge 0} t^{nk}.
\]
Hence, it follows that
\[
h_d(1,\zeta,\dots,\zeta^{n-1})=
  \begin{cases}
       1, & \text{if } n \mid d \\
       0, & \text{if } n\nmid d.\\
  \end{cases}
\]
The result is a consequence of the homogeneity property, since \[h_d(\lambda,\lambda\zeta,\dots,\lambda\zeta^{n-1})=\lambda^d \cdot h_d(1,\zeta,\dots,\zeta^{n-1}).\qedhere  \]  
\end{proof}

In this work, we use the Marar-Mond criterion \cite{Marar1989} to study some classes of corank one germs $f:(\CC^n,0)\rightarrow (\CC^p,0)$.

\begin{theorem}[{\cite[Theorem 2.14]{Marar1989}}]\label{thm: Marar-Mond crit}
For a finite germ $f:(\CC^n,0)\rightarrow (\CC^p,0)$ of corank $1$ 
 with $n < p$:
\begin{enumerate}[label=(\roman*)]
	\item $f$ is stable if and only if $D^k(f)$ is smooth of dimension $p-k(p-n)$, or empty, for any $k\geq1$.\label{item:mararmondi}
	\item \label{it2: Marar-Mond crit}$\eqA_e\textnormal{-codim}(f)$ is finite if and only if 
   for each $k$ with $p-k(p-n)\geq 0$, $D^k(f)$ is empty or an \textsc{icis} of dimension $p-k(p-n)$ and if, furthermore, for those $k$ such that $p-k(p-n)< 0$, $D^k(f)$ is $\{(0,\dots,0)\}$ or empty.\label{item:mararmondii}
\end{enumerate}
\end{theorem}

\section{Examples}\label{sec:Examples}
In this section, we present some examples to illustrate stable map germs. We will verify this property by calculating the $\mathcal{A}_e$-codimension of the germ or using the Marar–Mond criterion.
\begin{example}\label{ex:examplestable} The map germ $f:(\mathbb{C}^4,0) \to (\mathbb{C}^5,0)$ given by
\[f(x_1,x_2,x_3,y)=(x_1,x_2,x_3,y^2,x_1y)\]
is $\Ascr$-finite. Moreover, $f$ is stable.

In $f^*(\mathcal{O}_5)$ we obtain the following elements: $x_1^{a_1}x_2^{a_2}x_3^{a_3}y^{2k}$ and $ x_1^{a_1+1}x_2^{a_2}x_3^{a_3}y^{2k+1}$, for any non-negative integers $a_1,a_2$ and $k$. The extended tangent space $T\Ascr_e(f)$ has elements of the form 
\[\left(\begin{array}{cccc}
1 & 0 & 0 & 0\\
0 & 1 & 0 & 0\\
0 & 0 & 1 & 0\\
0 & 0 & 0 & 2y\\
y & 0 & 0 & x_1\\
\end{array}\right) \left(\begin{array}{c}
\epsilon_1\\
\epsilon_2\\
\epsilon_3\\
\epsilon_4\\
\end{array}\right)+\left(\begin{array}{c}
\eta_1(f)\\
\eta_2(f)\\
\eta_3(f)\\
\eta_4(f)\\
\eta_5(f)
\end{array}\right),\]
for any functions $\epsilon_i\in\Ocal_4$ ($i=1,\dots,4$) and $\eta_j\in\Ocal_5$ ($j=1,\dots,5$).
Also, for any non-negative integers $a_1,a_2$ and $k$, the following elements can be obtained in $T\Ascr_e(f)$:
\[\begin{aligned}
x_1^{a_1}x_2^{a_2}x_3^{a_3}y^{2k+1}e_1 &;\text{ by setting  } \left\{\begin{array}{l}
\epsilon_1=x_1^{a_1}x_2^{a_2}x_3^{a_3}y^{2k+1}, \epsilon_j =0 \text{ for } j\neq 1,\\[0.1cm]
\eta_5(f)=-x_1^{a_1}x_2^{a_2}x_3^{a_3}y^{2k+2}, \eta_l =0 \text{ for } l\neq 5,\vspace{0.2cm}\\
\end{array}\right. \\
x_1^{a_1}x_2^{a_2}x_3^{a_3}y^{2k+1}e_2
&;\text{ by setting} \ \epsilon_2=x_1^{a_1}x_2^{a_2}x_3^{a_3}y^{2k+1}, \epsilon_j=0 \text{ for } j\neq 2,\\[0.2cm]
x_1^{a_1}x_2^{a_2}x_3^{a_3}y^{2k+1}e_3
&;\text{ by setting  }\ 
\epsilon_3=x_1^{a_1}x_2^{a_2}x_3^{a_3}y^{2k+1}, \epsilon_j=0 \text{ for } j\neq 3,\\[0.2cm]
x_1^{a_1}x_2^{a_2}x_3^{a_3}y^{2k+1}e_4
&;\text{ by setting  } \left\{\begin{array}{l}
\epsilon_4=\frac{1}{2}x_1^{a_1}x_2^{a_2}x_3^{a_3}y^{2k},\epsilon_j=0 \text{ for } j\neq 4\\[0.1cm]
\eta_5(f)=-\frac{1}{2}x_1^{a_1+1}x_2^{a_2}x_3^{a_3}y^{2k}, \eta_l =0 \text{ for } l\neq 5, \vspace{0.2cm}\\
\end{array}\right. \\
x_1^{a_1}x_2^{a_2}x_3^{a_3}y^{2k+1}e_5
&;\text{ by setting  } \left\{\begin{array}{l}
\epsilon_1=x_1^{a_1}x_2^{a_2}x_3^{a_3}y^{2k} , \epsilon_j=0 \text{ for } j\neq 1,\\ [0.1cm]
\eta_1(f)=-x_1^{a_1}x_2^{a_2}x_3^{a_3}y^{2k}, \eta_l =0 \text{ for } l\neq 1.
\end{array}\right. \end{aligned}\]


From the previous comments it is clear that $T\Ascr_e(f)=\theta(f)$, so $f$ is stable. 
\end{example}

\quad

The example above is a particular case of the following general class of germs.

\begin{example}\label{ex:examplestableDk} The map germ $f:(\mathbb{C}^n,0) \to (\mathbb{C}^p,0)$ for $n < p < 2n$ given by
\[f(\xbm,y)=(x_1,\dots,x_{n-1},y^2,x_1y,\dots,x_{p-n}y)\]
is $\Ascr$-finite and, even more, $f$ is $\Ascr$-stable. We can check this property by the Marar-Mond criterion \cref{thm: Marar-Mond crit}.

Denoting the components of $f$ as in (\ref{Corank-1-germ}), the space $D^2(f)$ is given by the following divided differences:
\begin{itemize}
\item For $j=1$, the divided difference is
\[\begin{aligned}
f^{(2)}_{1}(\xbm,y_1,y_2)\coloneqq&\frac{f_{1}(\xbm,y_1)-f_{1}(\xbm,y_2)}{y_{1}-y_{2}}\\
=&\frac{y_1^2-y_2^2}{y_1-y_2}=y_1+y_2.
\end{aligned}\]
\item For $j=2,\dots,p-n+1$, the divided differences are
\[\begin{aligned}
f^{(2)}_j(\xbm,y_1,y_2)&\coloneqq\frac{f_j(\xbm,y_1)-f_j(\xbm,y_2)}{y_{1}-y_{2}}=\frac{x_{j-1}y_1-x_{j-1}y_2}{y_1-y_2}=x_{j-1}.
\end{aligned}\]
\end{itemize}

The space $D^3(f)$ is given by all the previous divided differences, $f_j^{(2)}$ for $j=1,\dots,p-n+1$, and the following new divided differences:
\begin{itemize}
\item For $j=1$, the new divided difference is
\[\begin{aligned}
f^{(3)}_{1}(\xbm,y_1,y_2,y_3)\coloneqq&\frac{f_{1}^{(2)}(\xbm,y_1,y_2)-f_{1}^{(2)}(\xbm,y_1,y_3)}{y_{2}-y_{3}}\\
=&\frac{y_1+y_2-y_1-y_3}{y_2-y_3}=1.
\end{aligned}\]

\item For $j=2,\dots,p-n+1$, the new divided differences are
\[\begin{aligned}
f^{(3)}_j(\xbm,y_1,y_2,y_3)&\coloneqq\frac{f_j^{(2)}(\xbm,y_1,y_2)-f_j^{(2)}(\xbm,y_1,y_3)}{y_{2}-y_{3}}=\frac{x_{j-1}-x_{j-1}}{y_2-y_3}=0.
\end{aligned}\]
\end{itemize}

Since $f^{(3)}_{1}$ is a unit, $D^k(f)$ is empty for all $k>2$. Thus, the double point scheme $$D^2(f)=\{(\xbm,y_1,y_2) \in \mathbb{C}^{n+1}; x_1=\ldots=x_{p-n}=0 \ \text{and} \ y_2=-y_1\}$$ is a smooth ICIS of dimension $2n-p$. Therefore, the Marar-Mond criterion guarantees that $f$ is stable.
\end{example}

\begin{example}
For the corank $2$ map germ $f(\xbm)=(x_1,x_2,y_1^2,y_2^2,y_1x_1,y_2x_1)$, where $\xbm=(x_1,x_2,y_1,y_2)$, the differential is
$$
df=\left(\begin{array}{cccc}
1 & 0 & 0 & 0\\
0 & 1 & 0 & 0\\
0 & 0 & 2y_1 & 0\\
0 & 0 & 0 & 2y_2\\
y_1 & 0 & x_1 & 0\\
y_2 & 0 & 0 & x_1
\end{array}\right).
 $$
It is clear that for any $k \geq 0$, then
\[y_2^{2k+1} e_3=\left(\begin{array}{c}
0\\
0 \\
y_2^{2k+1} \\
0 \\
0 \\
0
\end{array}\right)\text{ and }\ 
y_1^{2k+1} e_4=\left(\begin{array}{c}
0\\
0 \\
0 \\
y_1^{2k+1} \\
0 \\
0
\end{array}\right)
\]
are linearly independent in $N\Ascr_e(f)$, so $\codim_{\Ascr_e}(f)=\infty$. 
\end{example}

\section{Proof of Theorem A: the corank one case}\label{sec:corankone}

We present two proofs of \hyperlink{thm:monomialAfinite}{Theorem A} for corank one map germs. The first employs a geometric approach via the Marar-Mond criterion, and the second is purely algebraic and establishes the finiteness of the $\mathcal{A}_e$-codimension, providing a framework that extends to higher corank germs.

\quad
\hyperlink{thm:monomialAfinite}{Theorem A} states that a singular monomial map germ $f:(\CC^n,0)\to(\CC^p,0)$, with $n<p<2n$, is $\mathcal{A}$-finite if and only if $f$ is $\mathcal{A}$-equivalent to $$g(\xbm,y)=(x_1,\dots,x_{n-1}, y^2, x_{1}y,\dots,x_{p-n}y).$$
   
 In \cref{ex:examplestableDk}, we prove that $g$ is an $\Ascr$-stable map germ, hence $\Ascr$-finite (see also \cref{ex:examplestable}).

Now, consider a general corank one monomial map germ $f:(\CC^n,0)\to(\CC^p,0)$, with $n<p<2n$, given by 
\[ f(\xbm,y)=\big(x_1,\dots,x_{n-1},m_n(\xbm,y),\dots,m_{p}(\xbm,y)\big).\]

Our goal is to determine conditions on the monomials $m_i$ for $f$ to be $\Ascr$-finite.

\subsection{Proof using multiple point spaces}

Let $f:(\CC^n,0)\to(\CC^p,0)$ be a singular, $\Ascr$-finite monomial corank one map germ, with $n<p<2n$. In particular, $f$ is finite, that is, $f^{-1}(0)=\{0\}$. Then, by a change of coordinates, $f$ is given by 
\begin{equation}\label{eq:germform} f(\xbm,y)=\big(x_1,\ldots,x_{n-1},y^a,m_{n+1}(\xbm,y),\ldots, m_{p}(\xbm,y)\big),\end{equation}
for some $a \geq 2$.

We provide now a preliminary lemma that shows the importance of the instability hypothesis.

\begin{lemma}\label{prop: D2singular}
The set $D^2(f)$ is singular provided $f$ is unstable.
\end{lemma}
\begin{proof}
On the one hand, we know from \cref{thm: Marar-Mond crit} that some $D^k(f)$ is singular (or that some $D^k(f)\neq\varnothing$ for $p-k(p-n)<0$). 

On the other hand, in order to have a smooth $D^2(f)$ the set of divided differences must have linear terms that, if they are arranged as a vector, they have maximum rank (which is a simple reformulation of having differential of maximum rank).

Denoting by $m_n(\xbm,y)=y^a$ and $m_j(\xbm,y)=x_1^{a_{j,1}}\cdots x_{n-1}^{a_{j,n-1}}y^{b_j}$ for $j=n+1, \ldots, p$ and as shown in \cref{eq:bricks}, these divided differences have the form $$m_n^{(2)}(\xbm,y_1,y_2)=h_{a-1}(y_1,y_2) \ \ \text{and} \ \ m_j^{(2)}(\xbm,y_1,y_2)=x_1^{a_{j,1}}\cdots x_{n-1}^{a_{j,n-1}}h_{b_j-1}(y_1,y_2).$$ 
If they provide linear terms of maximum rank, the  divided differences for $D^3(f)$ must be either zero (if there is no variable $y_i$) or one (if there is): in the first case we do not have a space $D^3(f)$ of the correct dimension and in the other $D^3(f)$ would be empty. We have a contradiction in both cases, and the result follows by the Marar-Mond criterion \cref{thm: Marar-Mond crit}.
\end{proof}





To prove \hyperlink{thm:monomialAfinite}{Theorem A}, we separate the case $p=2n-1$ and $p<2n-1$.
\newline

\noindent
\textbf{Case $p=2n-1$.}

If $p=2n-1$, $D^2(f)\subset(\CC^{n+1},0)$ is an ICIS of dimension $p-2(p-n)=1$ and it is defined by $p-n+1=n$ equations. It may be non-irreducible, but by looking at the divided differences in \cref{eq:bricks}, each variable $x_i$ must appear in exactly one monomial. Indeed, one equation is of the form
\begin{equation*}
\sum_{a_1+a_2=a-1}y_1^{a_1}y_2^{a_2}\end{equation*}
and each of the remaining $n-1$ equations is of the form
\begin{equation*}
x_1^{\beta_1}\cdots x_{n-1}^{\beta_{n-1}}\sum_{a_1+a_2=a-1}y_1^{a_1}y_2^{a_2}.\end{equation*}
Therefore, if a variable appears in more than one monomial, we would have an irreducible component that does not have a variable in its defining equations and, therefore, be of dimension at least two: each of the irreducible components are defined by equations of the type
\[x_i^{\beta_i}=0\quad\text{and}\quad \sum_{a_1+a_2=a-1}y_1^{a_1}y_2^{a_2}.\] For the same reason, if a monomial is equal to $x_i^{\beta_i}y^{\alpha_i}$ with $\alpha_i>1$, the corresponding divided difference has a factor without the variable $x_i$ (see \cref{rem:factorization}), and we reach the same contradiction.

Hence, we have a map $\Ascr$-equivalent to
\[f(\xbm,y)=(x_1,\dots,x_{n-1},y^a,x_1^{\beta_1}y,\dots,x_{n-1}^{\beta_{n-1}}y).\]
The ICIS $D^2(f)$ is given by the equations
\[\begin{aligned}\sum_{a_1+a_2=a-1}y_1^{a_1}y_2^{a_2}&=0,\\
x_1^{\beta_1}&=0,\\&\vdots\\
x_{n-1}^{\beta_{n-1}}&=0.\end{aligned}\]
In turn, $D^3(f)$ is given by these previous equations and the new equations:
\begin{equation}\label{eq:d3divdiff}\begin{aligned}\sum_{a_1+a_2+a_3=a-2}y_1^{a_1}y_2^{a_2}y_3^{a_3}&=0,\\
0&=0,\\&\vdots\\
0&=0.\end{aligned}\end{equation}
If $a>2$ then $D^3(f)\subset(\CC^{n-1+3},0)$ does not have the correct dimension:  it has at least the same dimension as $D^2(f)\subset(\CC^{n-1+2},0)$ since its ambient space is one dimension higher but we only add one equation. Hence, $a=1,2$ and $D^3(f)$ is empty because for $a=2$ the only non-trivial equation would be $1=0$.

Moreover, $D^2(f)$ has a non-isolated singularity if some $\beta_i>1$. In conclusion, $\beta_i=1$ for all $i$, in which case $D^2(f)$ is smooth (or empty) and $D^3(f)$ is empty, which shows that $f$ is not unstable by  \cref{thm: Marar-Mond crit} and, even more, it must be stable of one of the forms given in \hyperlink{thm:monomialAfinite}{Theorem A}.
\newline

\noindent
\textbf{Case $p<2n-1$.}

 We suppose that $f$ of the form in \cref{eq:germform} is unstable and reach a contradiction. 
\newline

By \cref{prop: D2singular}, $D^2(f)$ must be a singular ICIS in variables $$x_1,\dots,x_{n-1},y_1 \ \text{ and } \ y_2$$ (with $p-n+1$ equations) and of dimension $p-2(p-n)\geq0$. If $y$ does not divide a coordinate function $m_i(\xbm,y)$, the corresponding divided difference is zero, which implies that the dimension of $D^2(f)$ cannot be $p-2(p-n)=n+1-(p-n+1)$. Furthermore, every variable must appear in some defining equation. If a variable does not appear in any equation, then $D^2(f)$ would be isomorphic to a product $V\times\CC$, where $V=\pi(D^2(f))$ and $\pi$ is the projection that removes the corresponding coordinate. This is incompatible with having isolated singularity. 

If $p<2n-1$, then $D^2(f)$ has dimension $p-2(p-n)=2n-p>1$ and since it is an ICIS, it is normal, hence reduced and irreducible. By \cref{rem:factorization}, irreducibility implies that every monomial must be either $m_i(\xbm,y)=x_j^{\beta_j}y$, for some $j$, or $m_i(\xbm,y)=y^2$ because a decomposition of one defining equation would provide a decomposition into factors of the ICIS. We reach a contradiction now because, as shown above, each of the variables must appear in some equation, but we only have $p-n+1<n$ divided differences (including that from $y^2$) of the form
\[x_j^{\beta_j}\quad \text{ or }\quad y_1+y_2.\]

In summary, we have shown that any corank one monomial map that is $\Ascr$-finite is stable. The normal form of \hyperlink{thm:monomialAfinite}{Theorem A} is a consequence of Marar-Mond criterion \cref{thm: Marar-Mond crit}, since $D^2(f)$ must be empty (and we have an embedding) or smooth of dimension $p-2(p-n)$ (and we have the other singular normal form). This is a straightforward computation.
\qed

\subsection{Proof using the $\Ascr_e$-codimension}
\label{Subsec:Ae-proof} We now present an alternative proof of Theorem A in the corank $1$ case, using Mather's infinitesimal criterion of finite determinacy, which asserts that a germ is $\Ascr$-finite if and only if its $\Ascr_e$-codimension is finite. The steps of this argument are essential for extending the result to the arbitrary corank case.

Let $f:(\CC^n,0)\to(\CC^p,0)$ be an $\Ascr$-finite corank one map germ, with $n<p<2n$. As before, by a change of coordinates, $f$ is given by 
\[ f(\xbm,y)=\big(x_1,\ldots,x_{n-1},y^a,m_{n+1}(\xbm,y),\ldots, m_{p}(\xbm,y)\big),\]
for some $a \geq 2$.

Let us analyze the following possibilities for the components $m_{j}$, $j\in \{n+1,\ldots,p\}$:
\newline

\noindent
\textbf{{1.-} Let $m_{n+1}(\xbm,y)=y^b$ with $gcd(a,b)=1$.}
   
In this case,
\[ f(\xbm,y)=\big(x_1,\ldots,x_{n-1},y^a,y^b,m_{n+2}(\xbm,y),\ldots, m_{p}(\xbm,y)\big).\]

We may suppose that $2 \leq a< b$, and that no coordinate is equal to a power of $y$ (in which case the argument is very similar). In this case, we can find any monomial $x_1^{a_1}\cdots x_{n-1}^{a_{n-1}}y^{\lambda a + \mu b}$ in $ f^*(\mathcal{O}_{p})$ for any non-negative integers $a_i,\lambda$ and $\mu$  ($i=1, \ldots, n-1$). 
       
Since $b \geq 3$, notice that $x_i^{d_i}ye_{n+1}$ cannot be reduced to zero in $  tf(\theta_n)$ under any relation in the quotient, for any $d_i \in \mathbb{N}$, and $i \in \{1, \ldots, n-1\}$. In order to ensure that $\codim_{\Ascr_e}(f)$ is finite, we must generate elements of the form
$$x_1^{\lambda_1}y, \ldots, x_{n-1}^{\lambda_{n-1}}y \in f^*(\mathcal{O}_{p})$$
for some $\lambda_i \in \mathbb{N}$. This requires adding at least $n-1$ distinct monomial components to $f$. However, we have only $p-n-1<2n-n-1=n-1$ remaining coordinate components available in $f$. Consequently, $\operatorname{codim}_{\mathcal{A}_e}(f) = \infty$, which contradicts the assumption that $f$ is $\mathcal{A}$-finite. Therefore, this case cannot occur.
\newline

\noindent
\textbf{2.- Let $m_{n+1}(\xbm,y)=x_1^{c}y^d$ with $c,d \in \mathbb{N}$.} 

Thus
\[ f(\xbm,y)=\big(x_1,\ldots,x_{n-1},y^a,x_1^cy^d,m_{n+2}(\xbm,y),\ldots, m_{p}(\xbm,y)\big).\]

Note that after renaming the variables and an $\Ascr$-equivalence, any variable $x_i$ could substitute for $x_1$ in $f_{n+1}$. We split this into two sub-cases:
\newline

\noindent
\textbf{2.1.-} If $a> 2$, then $x_j^{c_j}ye_{n} $ cannot be obtained from $ tf(\theta_n)$ (nor through any relation in the quotient)  for any $c_j \in \mathbb{N}$ and $j=2,\ldots, n-1$.

Hence, similarly to the previous case,
we must obtain 
$$x_2^{\lambda_2}y, \ldots, x_{n-1}^{\lambda_{n-1}}y \in f^*(\mathcal{O}_{p})$$
for some $\lambda_i \in \mathbb{N}$; so we must choose these elements as components of $f$. This forces $p=2n-1$ and the germ to be $\Ascr$-equivalent to
\[ f(\xbm,y)=\big(x_1,\ldots,x_{n-1},y^a,x_1^cy^d,x_2^{\lambda_2}y,\ldots, x_{n-1}^{\lambda_{n-1}}y\big).\]

 If $d>1$, then $x_1^{c_1}y e_{n} \notin T\Ascr_e(f)$ for any $c_1 \in \mathbb{N}$. If $d=1$ and $c>1$ then, in particular, $y^b e_{n+1} \notin T\Ascr_e(f)$ for any $b \not\equiv 0 \pmod{a}$. 
 
 We still need to analyze the case $c=d=1$ and, by symmetry, $\lambda_i=1$; i.e., the previous germ must be $\Ascr$-equivalent to \begin{equation}\label{eq:refhere} f(\xbm,y)=(x_1,\ldots,x_{n-1},y^a, x_1y,\ldots,x_{n-1}y)\end{equation} 
 with $a>2$. 

To ensure that all but finitely many powers $y^{\delta}e_{n+1}$ belong to $T\Ascr_e(f)$  we must use $tf$, since for
$\delta_1 \not\equiv 0 \pmod{a}$ we cannot obtain it only from $f^*(\Ocal_p)$. The only vector field that obtains $y^{\delta}e_{n+1}$ as a component by $tf$ (see \cref{Tang-Space}) is 
\[\epsilon=y^{\delta -1}e_1  \ \mapsto \ tf(\epsilon)= y^{\delta -1}e_1 + y^\delta e_{n+1}.\]

However, we cannot eliminate the term $y^{\delta_1-1}$ using $tf$ again nor by an element in $f^*(\mathcal{O}_{p})$, for $\delta_1 \equiv 2 \not\equiv 0 \pmod{a}$. In this case, $f$ would not be $\Ascr$-finite.
\newline

 \noindent
\textbf{2.2.-} Let $a=2$, so the germ has the form  \[ f(\xbm,y)=\big(x_1,\ldots,x_{n-1},y^2,x_1^cy^d,m_{n+2}(\xbm,y),\ldots, m_{p}(\xbm,y)\big).\]  
If $d>1$, we cannot obtain $x_j^{\delta_j}ye_{n+1}$ in $tf(\theta_n)$ for any $\delta_j \in \mathbb{N}$ with $j =2,\dots,n-1$. Nor can these elements be obtained in $f^*(\mathcal{O}_p)$, unless they appear as components of $f$, that is, $p=2n-1$ and \[(\xbm,y)=\big(x_1,\ldots,x_{n-1},y^2,x_1^cy^d,x_2^{\lambda_2}y,\ldots, x_{n-1}^{\lambda_{n-1}}y\big),\] for some $\lambda_2, \ldots, \lambda_{n-1} \in \mathbb{N}$. If $c>1$, the element $y^{2k+1}e_{n+1}$ does not belong to $T\mathcal{A}_e(f)$ for any $k \in \mathbb{N}$. Since $f$ has no remaining free components, these cases cannot occur. This argument yields $c=1$. Finally, if $d$ is even, the term $x_1y^d$ is eliminated;  if $d>1$ is odd, the element $x_1^{\delta_1}ye_{n+1}\notin T\mathcal{A}_e(f)$. Therefore, for $c>1$ or $d>1$, the map germ $f$ is not $\Ascr$-finite.
    
If $c=1=d$ and the coordinate functions $m_{i}$ ($i=n+2,\dots,p$) depend only on two variables, an argument similar to the one made in the previous sub-case (see \cref{eq:refhere}) shows that $f$ must be $\Ascr$-equivalent to
\[(\xbm,y)\mapsto(x_1,\ldots,x_{n-1}, y^2,  x_{1}y,\ldots,x_{p-n}y),\]
as desired.
     
If $d=1=c$ and $m_{n+2}$ is a monomial in three variables or more (say $x_i, x_j$ and $y$), then we cannot obtain  $y^{2\lambda +1}e_{n+2}$ in $tf(\theta_n)$ for any $\lambda \in \mathbb{N}$. We cannot obtain them from $f^*(\Ocal_p)$ either, since no other coordinate function $f_{n+i}(\xbm,y)$ can be equal to $y^b$, as shown in case 1. This shows the contradiction.
    \newline

\noindent 
\textbf{3.-} If $m_{n+1}$ depends on three variables or more, as we analyzed in the last subcase, $f$ cannot be $\Ascr$-finite: the element $y^{2\lambda +1}e_{n+1}$ can be obtained neither in $tf(\theta_n)$ for any $\lambda \in \mathbb{N}$ nor from $f^*(\Ocal_p)$.\qed

\section{Proof of Theorem A: General case}\label{sec:coranktwo}

We now show that no monomial map that is $\Ascr$-finite can have corank greater than one, thereby completing the proof of \hyperlink{thm:monomialAfinite}{Theorem A}. More precisely, we suppose that $f:(\CC^n,0)\to(\CC^p,0)$ is a monomial $\Ascr$-finite map with rank $0 \leq r<n-1$ and $n< p<2n$, and reach a contradiction.
\newline

Since $f$ is finite, it must be $\Ascr$-equivalent to
\[f(\xbm,\ybm)=\big(x_1,\dots,x_r,y_1^{a_1},\dots,y_{n-r}^{a_{n-r}},m_{n+1}(\xbm,\ybm),\dots,m_p(\xbm,\ybm)\big),\]
where $\xbm = (x_1,\ldots, x_r)$ and $\ybm =(y_1,\ldots, y_{n-r})$. The corank of $f$ is greater than one, so $1< n-r \leq n$. Up to a relabeling of the variables, if necessary, we suppose $2 \leq a_1 \leq \ldots \leq a_{n-r}$. 

Finally, abusing our notation, we will say that two (or more vectors) are linearly independent in $N\Ascr_e(f)$ if their classes are so after taking the quotient by $T\Ascr_e(f)$.
\newline

\noindent
\textbf{Step 1: Consider $a_1\geq3$.}

With a similar argument as in Subsection \ref{Subsec:Ae-proof}, item 1, we show that the sets 
\[\big\{x_i^{\beta_i} y_1 e_{r+1}:\beta_i\in\Nbb, i=1,\dots,r\big\},\] 
\[\quad\big\{y_j^{\beta_j} y_1 e_{r+1}:\beta_j\in\Nbb, j=2,\ldots, n-r\big\},\  \text{and}\]
\[
\big\{ y_2^{a_2\beta+1} e_{r+1}:\beta\in\Nbb\big\}\]
contain infinitely many linearly independent classes in $N\Ascr_e(f)$.

Since $a_1\geq 3$, none of the elements from the previous sets belong to $tf(\theta_n)$ via any relation in the quotient. To ensure that $\operatorname{codim}_{\mathcal{A}_e}(f)$ is finite, we must include at least $n$ distinct monomials from these sets
\[\big\{x_i^{\beta_i} y_1:\beta_i\in\Nbb,i=1,\dots,r\big\}\cup\big\{y_j^{\delta_j} y_1 :\delta_j\in\Nbb, j=2,\dots,n-r\big\}\cup\big\{ y_2^{a_2\kappa+1} :\kappa\in\Nbb\big\}\]
as coordinate functions of $f$. For instance:

\begingroup
\crefformat{enumi}{(#2#1#3)}
\Crefformat{enumi}{(#2#1#3)}
\begin{enumerate}
    \renewcommand\theenumi{(\alph{enumi})}
    \renewcommand\labelenumi{\theenumi}
    \item \label{en:type1}$x_i^{b_i}y_1$, for some natural numbers $b_i$ and all $i=1,\dots,r$;
    \item \label{en:type2}$y_j^{d_j}y_1$ for some natural numbers $d_j$ with $j=2,\dots,n-r$;
    \item $y_2^{a_2 k +1}$, for some natural number $k$.
\end{enumerate}
\endgroup

However, $f$ has only $p-n<n$ remaining free components. This conclusion contradicts our assumption that $f$ is $\Ascr$-finite.
\newline

By a symmetry argument, $a_i=2$ for all $i$, and $f$ must be $\Ascr$-equivalent to
\[f(\xbm,\ybm)=\big(x_1,\dots,x_r,y_1^2,\dots,y_{n-r}^2,m_{n+1}(\xbm,\ybm),\dots,m_p(\xbm,\ybm)\big).\]

\vspace{0.2cm} \noindent
\textbf{Step 2: Consider $m_{n+1}(\xbm,\ybm)= y_l^N$ for $N>1$ and some $l \in \{1,\ldots, n-r\}$.}

Up to an $\Ascr$-equivalence, we can assume that $m_{n+1}(\xbm,\ybm)=y_1^N$, for some $N \not\equiv 0 \pmod{2}$. 

There is no way to obtain any vector $h(\xbm,y_2,\dots,y_{n-r}) y_1 e_{n+1}$ in $T\Ascr_e(f)$ from the differential part of the quotient, that is, as a term in  $df(\omega)$ for some $\omega\in\theta_n$, since $N\geq3$. Hence, they must be obtainable from $f^*(\Ocal_n)$, and some coordinate functions of $f$ must be equal to those shown in \cref{en:type1,en:type2} above. This, again, exceeds the number of coordinates available in $f$.

\vspace{0.2cm}\noindent
\textbf{Step 3: Consider $m_{n+1}(\xbm,\ybm)=y_1^c y_2^d$, for some $c,d\in \mathbb{N}$.}

Naturally, the same argument holds for $m_{n+1}(\mathbf{x},\mathbf{y})=y_i^c y_j^d$, for some $c,d\in \mathbb{N}$ and $i,j\in \{1,\ldots, n-r\}$ with $i\neq j$.

If $d > 1$, we cannot obtain from $tf(\theta_n)$, via any relation in the quotient, the following $n-2$ elements : $x_i^{\delta_i}y_1e_{n+1}$ and $y_j^{\gamma_j}y_2e_{n+1}$, for any $\delta_i, \gamma_j \in \mathbb{N}$ with $i = 1, \dots, r$ and $j=3,\ldots, n-r$.  

Nor can these elements be obtained from $f^*(\mathcal{O}_p)$, unless they appear as components of $f$, for instance: 
\begingroup
\crefformat{enumi}{(#2#1#3)}
\Crefformat{enumi}{(#2#1#3)}
\begin{enumerate}
    \renewcommand\theenumi{(\alph{enumi})}
    \renewcommand\labelenumi{\theenumi}
    \item[(a1)] \label{en:type-a1}$y_j^{\mu_j}y_2$ for some natural numbers $\mu_j$ with $j=3,\dots,n-r$;
    \item[(a2)] \label{en:type-a2}$x_i^{\lambda_i}y_1$, for some natural numbers $\lambda_i$ and all $i=1,\dots,r$.
\end{enumerate}
\endgroup
Thus, $p=2n-1$ and 
$$f(\xbm,\ybm)=\big(x_1,\ldots,x_r,y_1^2,\ldots, y_{n-r}^2,y_1^cy_2^d,y_3^{\mu_3}y_2,\ldots, y_{n-r}^{\mu_{n-r}}y_2,x_1^{\lambda_1}y_1,\ldots, x_r^{\lambda_r}y_1\big).$$

For $d>1$, if 
\[\left\{\begin{array}{ll}
   c>1, & \mbox{then}\ y_1^{2k+1}e_{n+1}, y_2^{2r+1}e_{n+1} \notin T\mathcal{A}_e(f) \medskip\\
   c=1, & \mbox{then}\ y_1^{2k+1}e_{n+1} \notin T\mathcal{A}_e(f), 
\end{array} \right.\]
for any relation in the quotient and for any $k,r \in \mathbb{N}$. Then, $f$ would not be $\Ascr$-finite, since it has no remaining free components. 

Thus $c=d=1$ and, by the same reasoning as before, $\mu_3=\ldots=\mu_{n-r}=1$. However, $y_2^{2k+1}e_{r+1}$ cannot be generated in $T\mathcal{A}_e(f)$ for any $k \in \mathbb{N}$.

Therefore, in every case considered in Step 2, the germ $f$ fails to be $\Ascr$-finite. 

An analogous argument applies if $m(\mathbf{x},\mathbf{y})=x_i^cy_j^d$.

\vspace{0.2cm}\noindent
\textbf{Step 4: If $m_{n+1}$ is a monomial in three or more variables.}

As analyzed in the previous step, $f$ is not $\mathcal{A}$-finite, since we do not have enough free components to accommodate the monomials $q(\xbm,\ybm)$ whose vector $q(\xbm,\ybm)e_{n+1}$ was obtained neither in $tf(\theta_n)$ nor in $f^*(\Ocal_p)$.\qed
\newline

An immediate corollary of \hyperlink{thm:monomialAfinite}{Theorem A} is that a singular monomial parameterized hypersurface in $\CC^{n+1}$ is $\Ascr$-finite if and only if it is $\Ascr$-equivalent to a trivial unfolding of the cross-cap.

\begin{corollary}
Let $f:(\CC^n,0)\to(\CC^{n+1},0)$ be a monomial map germ. Then $f$ is $\Ascr$-finite if and only if $f$ is $\Ascr$-equivalent to one of the following $\Ascr$-stable map germs:
\begin{enumerate}
    \item $(\xbm,y)\mapsto(x_1,\ldots, x_{n-1},y,0)$,
    \item $(\xbm,y)\mapsto(x_1,\ldots,x_{n-1}, y^2,  x_{1}y)$.
\end{enumerate} 
\end{corollary}

\section{Proof of Theorem B}\label{sec:ngep}

In this section we prove \hyperlink{thm:ngeqp}{Theorem B}. Together with \hyperlink{thm:monomialAfinite}{Theorem A}, this completes the classification of $\Ascr$-finite monomial map germs $f:(\CC^n,0)\to(\CC^p,0)$ for all $p<2n$ (see \cref{cor:p<2n}).

\subsection{Reduction and $\Kscr$-equivalence}

We use the \textit{Mather-Gaffney geometric criterion}, which shows a geometric condition that is equivalent to being $\Ascr$-finite.

\begin{theorem}[see {\cite[Theorem~2.1]{Wall1981}} and {\cite[Theorem~4.5]{Mond2020}}]
A holomorphic map germ $f:(\CC^n,0)\to(\CC^p,0)$ has finite $\Ascr_e$-codimension if and only if there is a small enough representative $f:X\to Y$ such that:
\begin{enumerate}
\item $f^{-1}(0)\cap C(f)\subset\{0\}$;
\item the restriction $f:X\setminus f^{-1}(0)\to Y\setminus\{0\}$ is locally stable;
\end{enumerate}
where $C(f)$ denotes the critical set of $f$.
\end{theorem}

Moreover, it is a standard result that the first condition, $f^{-1}(0)\cap C(f)\subset\{0\}$, is equivalent to $\Kscr$-finiteness. This is a consequence of the coherence of the sheafification of the $\Kscr_e$-tangent space and of \cite[Theorem E.3]{Mond2020} (see  \cite[Exercise 4.5.9]{Mond2020}). Therefore, the Mather-Gaffney criterion shows that $\Ascr$-finiteness implies $\Kscr$-finiteness.

\begin{remark}\label{rem:dimf0}
If $n\geq p$ and $f^{-1}(0)\cap C(f)\subset\{0\}$, then $f^{-1}(0)$ has dimension $n-p$. Indeed, since $f^{-1}(0)$ is defined by at most $p$ equations, $\dim f^{-1}(0) \geq n-p$. Let $Z$ be an irreducible component of $f^{-1}(0)$ and assume that $\dim Z>n-p\geq0$. Then $Z\setminus\{0\}\neq\varnothing$ and $f$ is not a submersion at any $z\in Z\setminus\{0\}$ (otherwise, $f^{-1}(0)$ would be a manifold of dimension $n-p$ around $z$). Therefore, $Z\setminus\{0\}\subseteq f^{-1}(0)\cap C(f)$, a contradiction.
\end{remark}

What we do now is reduce all the possible $\Ascr$-finite monomial map germs by classifying all possible monomial $\Kscr$-finite germs.

\begin{lemma}\label{lem:Kfinitemonomial}
Let $f:(\CC^n,0)\to(\CC^p,0)$ be a monomial map germ with $n\geq p$. Then, $f$ is $\Kscr$-finite if and only if, after a permutation of the coordinates of the source and the target, $f$ is one of the following germs:
\begin{enumerate}[label=(\alph*)]
    \item\label{it:Ksubm} $(x_1,\dots,x_n)\mapsto(x_1,\dots,x_p)$;
    \item\label{it:Kpowers} $(x_1,\dots,x_n)\mapsto(x_1^{a_1},\dots,x_n^{a_n})$, with $a_i\geq1$  and $n=p$;
    \item\label{it:Kmorse} $(x_1,\dots,x_n)\mapsto(x_1,\dots,x_{p-1},x_px_{p+1})$, with  $n=p+1$.
\end{enumerate}
\end{lemma}
\begin{proof}
First, let $f$ be a monomial $\Kscr$-finite map germ, given as follows
\[f(\xbm)=\big(f_1(\xbm),\dots,f_p(\xbm)\big).\]
By \cref{rem:dimf0}, if more than one coordinate $f_i$ depends on the same variable $x_j$, then $f^{-1}(0)$ has an irreducible component contained in $V(x_j)$ with at most $p-1$ equations, so $\dim f^{-1}(0)>n-p$, a contradiction. This shows that the coordinate $f_i$ depends on some variables $\xbm_i=(x_{i,j})_j$, $i=1,\dots,p$, and they are all pair-wise different. Possibly, there are some other free variables $\xbm_F=(x_{F,j})_j$ that do not appear in any coordinate:
\[f(\xbm_1,\dots,\xbm_r,\xbm_F)=\big(f_1(\xbm_1),\dots,f_p(\xbm_r)\big).\]

Since all coordinates depend on different variables, the differential of $f$ has blocks $\nabla f_i(\xbm_i)$. Therefore, a critical point of $f$ is the union of all critical points of $f_i(\xbm_i)$.

If no coordinate $f_i$ is critical, then we are in case $(a)$.

If one coordinate is critical at the origin (say $f_1$), then all other coordinates can depend on at most one variable. Indeed, if a different coordinate (say $f_2$) depends on more than one variable, then all the points that have zero in all coordinates except $x_{2,1}$ (which is free) are in $f^{-1}(0)$ and also in $C(f)$. This would contradict $\Kscr$-finiteness. The same argument shows that there are no free variables $\xbm_F$:
\[f(\xbm_1,x_s\dots,x_{p},\xbm_F)=\big(f_1(\xbm_1),x_s^{a_s},\dots,x_p^{a_p}\big),\] for some $s$.
However, exactly the same argument shows that if some $a_i$ is greater than one, then $f_i$ is critical and $f_1$ must depend on one variable only, so we are in case $(b)$.

The only case left is that $f_1$ is critical and all $a_i=1$:
\[f(\xbm_1,x_s\dots,x_{p},\xbm_F)=\big(f_1(\xbm_1),x_s,\dots,x_p\big).\] 
If $f_1$ depends on more than 2 variables, at least $x_{1,1}$, $x_{1,2}$ and $x_{1,3}$, then $\{x_{1,1}=x_{1,2}=0\}$ is contained in the zero set of $f_1$ and also in the critical set of $f_1$, so $f$ cannot be $\Kscr$-finite. Similarly, if $f_1$ only depends on two variables and 
\[x_{1,1}^2x_{1,2}\]
divides $f_1$, then $\{x_{1,1}\}$ is also contained in $f_1^{-1}(0)\cap C(f_1)$, so $f$ cannot be $\Kscr$-finite either. The only possibility left is case $(c)$.
\newline

It is easy to check that these three cases satisfy $f^{-1}(0)\cap C(f)\subseteq\{0\}$, so they are all $\Kscr$-finite.
\end{proof}

\subsection{$\Ascr$-finiteness}

To prove \hyperlink{thm:ngeqp}{Theorem B}, we only need to see which cases from \cref{lem:Kfinitemonomial} are $\Ascr$-finite.
\newline

Case $(a)$ is obviously $\Ascr$-finite because it is stable.
\newline

For case $(b)$, by \cref{Tang-Space}, the tangent space $T\Ascr_e f$ is given as follows.
\begin{equation}
   T\Ascr_e(f)=\left\{\left(\begin{array}{c}
a_1\epsilon_1 x_1^{a_1-1}+\eta_1(x_1^{a_1},\dots,x_n^{a_n})\\
\vdots\\
a_n\epsilon_n x_n^{a_n-1}+\eta_p(x_1^{a_1},\dots,x_n^{a_n})
\end{array}\right); \begin{array}{c}
\epsilon \in \theta_n, \ \eta \in \theta_n,\\
\epsilon=(\epsilon_1,\ldots,\epsilon_n), \\
\eta=(\eta_1,\ldots,\eta_n)
\end{array}\right\}.
\end{equation}
It is easy to see that if more than one $a_i$ is bigger than one (say $a_{n-1},a_n>1$), then $x_{n-1}^{k a_{n-1}+1}e_n$ is not in $T\Ascr_e f$ for any $k$, so $f$ is not $\Ascr$-finite. Similarly, if $a_n>2$ and there are more than one variable ($n>1$), then $x_1^k x_ne_n$ is not in $T\Ascr_e f$ for any $k$, so $f$ is not $\Ascr$-finite.

This only leaves two cases, $f(x)=x^a$ and $n=p=1$, or
\[f(x_1,\dots,x_n)=(x_1,\dots,x_{n-1},x_n^2)\]
and $n=p>1$. It is easy to check that the former has $\Ascr_e$-codimension $a-2$ and the latter is stable (it is a fold-type stability).

\vspace{0.2cm}
It is an easy exercise to see that $(c)$ is stable.
\vspace{0.2cm}

This concludes the proof of \hyperlink{thm:ngeqp}{Theorem B}.
\newline

\begin{corollary}\label{cor:p<2n}
Let $f:(\CC^n,0)\to(\CC^p,0)$ be a monomial map germ with $p<2n$ and $(n,p)\neq(1,1)$. Then $f$ is $\Ascr$-finite if and only if it is $\Ascr$-stable.
\end{corollary}

\section{Non-monomial $\Ascr$-finite map germs}\label{sec:Non-monomial}

In this section, we study non-monomial germs. We begin by characterizing stable germs of corank 1. Subsequently, we construct families of non-monomial $\mathscr{A}$-finite germs associated with monomial ones.

\subsection{Stable singularities}

As seen in \hyperlink{thm:monomialAfinite}{Theorem A}, an $\Ascr$-finite monomial map germ $f : (\CC^n,0)\to(\CC^p,0)$ with $n < p < 2n$ is necessarily $\mathcal{A}$-stable, being either an immersion or a map germ of multiplicity two.

After the work of Mather in \cite{Mather1969a}, stable map germs have been classified by their local algebras. Although this characterization is widely used to construct stable germs, only \textit{ad hoc} methods have been available for the specific cases considered here. For instance, in \cite[p. 140]{Cooper2002}, \cite[Theorem 3.1]{Houston2005b}, and \cite[Proposition 3.4]{OsetSinha2022}, non-canonical forms were used to relate such germs to sections of $\mathcal{A}_e$-codimension one.

 In \cite{Morin1965}, Morin shows a result that, in modern language, implies that $f: (\CC^{n},0)\to (\CC^{n+1},0)$ is a stable corank one map germ if and only if there exists an integer $d$ with $2\leq d \leq \frac{n+2}{2}$ such that $f$ is $\Ascr$-equivalent to \[ \begin{array}{ccl}
f_d : (\CC^{n-2d+2}\times\CC^{d-2}\times\CC^{d-1}\times\CC,0) & \to &  (\CC^{n-2d+2}\times\CC^{d-2}\times\CC^{d-1}\times\CC^2,0)\\
 (\wbm,\ \ubm,\ \vbm, \ y) &\mapsto &   \left(\wbm,\ubm,\ \vbm,\ y^d+\displaystyle\sum_{i=1}^{d-2}u_iy^i, \ \sum_{i=1}^{d-1}v_iy^i \right),\end{array}\] with $\wbm \in \CC^{n-2d+2}$, $\ubm=(u_1,\ldots, u_{d-2})\in \CC^{d-2}$ and $\vbm=(v_1,\ldots,v_{d-1})\in \CC^{d-1}$.


\begin{theorem}\label{thm:corankonestables}
Let $f : (\mathbb{C}^n, 0) \to (\mathbb{C}^p, 0)$ be a corank one map germ, with $n < p \le 2n$. The germ $f$ is $\mathscr{A}$-stable if and only if there exists an integer $d \ge 1$ satisfying $(p - n + 1)(d - 1) \le n$ such that $f$ is $\mathscr{A}$-equivalent to the trivial unfolding $f_d$, given by
\[
f_d(\ubm, y) = \left(\ubm,\ y^d + \sum_{i=1}^{d-2} x_{1,i} y^i,\ \sum_{i=1}^{d-1} x_{2,i} y^i,\ \ldots,\ \sum_{i=1}^{d-1} x_{p-n+1,i} y^i \right),\]
where $\ubm=(\wbm, \xbm_1, \xbm_2, \ldots, \xbm_{p-n+1})\in \CC^{n-1}$, $\wbm \in \mathbb{C}^{n - (p-n+1)(d-1)}$, $\xbm_1 = (x_{1,1}, \ldots, x_{1,d-2}) \in \mathbb{C}^{d-2}$, $\xbm_j = (x_{j,1}, \ldots, x_{j,d-1}) \in \mathbb{C}^{d-1}$ for each $j = 2,\ldots, p-n+1$.
\end{theorem}
\begin{proof}
Since stable germs are $\Ascr$-equivalent if and only if they are $\Kscr$-equivalent (see \cite{Mather1969a}, cf. \cite[Chapter 7]{Mond2020}), we must find a stable germ whose local algebra isomorphic to 
\[\frac{\Ocal_1}{(y^{d})},\]
for some $d$, because we are assuming that the germ has corank one.  Since unfoldings of germs have the same local algebra as that of the germ they unfold, it is enough to find a \textit{minimal} stable germ with such a local algebra (i.e., a germ whose isosingular locus is just a point) because any other will be a trivial unfolding of a minimal one. Observe that the germ of the statement is a trivial unfolding with parameters $\wbm$.

Moreover, since multiple point spaces of stable germs are smooth complete intersections by the Marar-Mond criterion \cref{thm: Marar-Mond crit}, we must first see what the equations of such multiple point spaces are. By smoothness, the linear parts of the equations provide a matrix of maximal rank, so by a change of variables, an $\Rscr$-equivalence, we can assume that the equations are
\[\begin{aligned}
0&=x_i+h.o.t.,\text{ or }\\
0&=y_1+\dots+y_{s}+h.o.t.
\end{aligned}\]
This $\Rscr$-equivalence on the multiple point spaces will be induced by an $\Ascr$-equivalence on the germ (actually, another $\Rscr$-equivalence) if the $\Rscr$-equivalence on the divided differences is symmetric with respect to the variables $y_j$, which we can always assume in this case.

Observe that if one equation is of the form
\[0 =y_1+\dots+y_{s}+h.o.t.,\]
its divided difference will be a unit. This shows that there is only one such equation in the set of divided differences that defines the (nonempty) multiple point spaces. Thus, $s=d-1$ and the 
 minimal germ must be ($\Ascr$-equivalent to)
\[f(\xbm,y)=\left(\xbm,\ y^{d}+\sum_{i=1}^{d-2}x_{1,i}y^i,\ \sum_{i=1}^{d-1}x_{2,i}y^i,\ \dots,\ \sum_{i=1}^{d-1}x_{p-n+1,i}y^i\right),\]
 because it has the divided differences we required.

 Since we were considering only minimal stable germs, the germ $f_d$ of the statement is a trivial unfolding of $f$, with parameters $\wbm \in \mathbb{C}^{n - (p-n+1)(d-1)}$.
\end{proof}

\subsection{$\Ascr$-finite unfoldings}

One strategy to obtain non-monomial, $\Ascr$-finite surface map germs in $\mathbb{C}^3$ is to consider suitable deformations (unfoldings) of monomial map germs that satisfy the finite determinacy property.

More precisely, starting from a (not necessarily $\Ascr$-finite) monomial map germ
\[
g(y)=(y^a,y^b),
\]
we consider 1-parameter deformations of the form
\[
f(x,y)=\big(x,g_x(y)\big)
\]
that are $\Ascr$-finite.

While this is a significant challenge in full generality, we restrict our attention to study these deformations for $a$ equal to $2$ or $3$ and such that
\[g_x(y)=\big(y^a,y^b+m(x,y)\big)\]
for some monomial $m(x,y)$. 

This general framework recovers, for instance, the well-known classification by Mond \cite{Mond1985}. He proved that any simple map germ from $(\mathbb{C}^2,0)$ to $(\mathbb{C}^3,0)$ is $\Ascr$-equivalent to one in the list:
    $$
\begin{array}{rll}
          & f(x,y) = (x, y, 0) & \text{(immersion)} \\
    S_0 : & f(x,y) = (x,y^2, xy) & \text{(cross cap)} \\
    S_k : & f(x,y) = (x, y^2, y^3 + x^{k+1} y), & k \geq 1 \\
    B_k : & f(x,y) = (x, y^2, x^2 y + y^{2k+1}), & k \geq 2 \\
    C_k : & f(x,y) = (x, y^2, xy^3 + x^k y), & k \geq 3 \\
    F_4 : & f(x,y) = (x, y^2, x^3 y + y^5) & \\
    H_k : & f(x,y) = (x, y^3, xy + y^{3k-1}), & k \geq 2.
\end{array}
$$ A map germ singularity is called \textit{simple} if only finitely many $\Ascr$-classes of singularities can be found in all possible perturbations.

In the next result, except for the $C_k$ family and the immersion (which is trivial), we recover all simple map germs classified by D. Mond by deforming a plane curve germ $g(y) = (y^a, y^b)$ via a specific 1-parameter unfolding.

\begin{theorem}\label{prop:defofmonomial}
Let $m \in \mathbb{C}[x,y]$ be a monomial and $f:(\mathbb{C}^2,0)\to (\mathbb{C}^{3},0)$ be a singular map germ given by \[f(x,y)=\big(x, y^a, y^b+m(x,y)\big),\] with $a=2$ or $a=3$, and $b>a$. Then, $f$ is $\Ascr$-finite if and only if $f$ is $\Ascr$-equivalent to:  
\begin{enumerate}
\item[(i)] $(x,y)\mapsto\big(x, y^2, xy)$,
\item[(ii)] $(x,y)\mapsto\big(x, y^2, y^{2k+1}+x^\lambda y\big)$, or
\item[(iii)] $(x,y)\mapsto\big(x, y^3, y^{3k+2}+x^\lambda y\big)$,
\end{enumerate} for some positive integer $\lambda$.
\end{theorem}

\begin{proof}
Since $m$ is a monomial, let $f(x,y)=(x, y^a, y^b+\alpha x^{\lambda}y^{\mu})$ be an $\Ascr$-finite map germ, for some $\lambda,\mu \ge 0$, and $\alpha \in \mathbb{C}$. 

We can now assume, up to $\Ascr$-equivalence, that the monomial $m(x,y)$ has coefficient one, since, by \hyperlink{thm:monomialAfinite}{Theorem A}, the map germ $(x,y) \mapsto (x, y^a, y^b)$ is not $\Ascr$-finite.

We analyze each congruence class of $b$ modulo $a$, for $a=2$ or $a=3$.

For $b=ak$, it is easy to see that $f(x,y) = (x, y^a, y^{ak} + x^{\lambda}y^{\mu})$ is $\Ascr$-equivalent to $g(x,y) = (x, y^a, x^{\lambda}y^{\mu})$. 

 Since $g$ is a monomial germ, by \hyperlink{thm:monomialAfinite}{Theorem A}, $g$ is $\Ascr$-finite if and only if it is $\Ascr$-equivalent to the cross-cap, that is, $(x,y) \mapsto (x, y^2, xy)$.

For each $a \in \{2,3\}$ and $b=ak+r$ with $1\leq r< a$, we suppose that the  map germ  \[
f(x,y)=(x, y^a, y^{ak+r}+x^{\lambda}y^{\mu})
\] is $\Ascr$-finite. Note that $\mu \not\equiv 0 \pmod a$, otherwise $f$
is $\Ascr$-equivalent to
$(x,y) \mapsto (x, y^a, y^{ak+r})$, which is not $\Ascr$-finite, by the previous argument.

If $\mu >2$, then $f$ is not $\Ascr$-finite because $x^cye_3 \notin T\Ascr_e(f)$ for any $c\geq 0$ (similarly for any other relation in the quotient where it may appear). Likewise, if $\mu = 2$ and $a=3$ (as for $a=2$ the argument reduces to the previous case), then $x^cye_2 \notin T\Ascr_e(f)$ for any $c \geq 0$. In both cases the $\Ascr_e$-codimension is not finite, contradicting the hypothesis that $f$ is $\Ascr$-finite. Therefore, we conclude that $\mu =1$.

If $\lambda =0$, then $f(x,y)=(x, y^a, y^{ak+r}+y)$ is immersive. However, $f$ is a singular map germ. 

Thus, $f$ is $\Ascr$-equivalent to one of the following map germs:
\begin{enumerate}
    \item $(x,y) \mapsto (x,y^2,y^{2k+1}+x^{\lambda}y)$,
    \item $(x,y) \mapsto (x,y^3,y^{3k+1}+x^{\lambda}y)$,
    \item $(x,y) \mapsto (x,y^3,y^{3k+2}+x^{\lambda}y)$.
\end{enumerate}
However, germ (2) does not occur by \cref{prop-germ-3-3k+1}. 

Reciprocally, the $\mathcal{A}$-finiteness of the remaining germs in cases (1) and (3) is established in \cref{prop:mult2} and \cref{prop:mult3}, respectively.
\end{proof}

Indeed, the cases shown in \cref{prop:defofmonomial} above are the only $\Ascr$-finite map germs of multiplicity $2$ or $3$ that can be obtained in this manner (with the exception of the immersion). We show now the non-$\Ascr$-finiteness of the  germ (2) in the proof of the previous result.

 \begin{lemma}
 \label{prop-germ-3-3k+1}
The map germ $f(x,y)=(x,y^3,y^{3k+1}+x^\lambda y)$
is not $\Ascr$-finite.
 \end{lemma}
\begin{proof}
Let $h_d$ be the complete homogeneous symmetric polynomial of degree $d$, as in \cref{def:Homog-Polyn-d}. The iterated divided differences of this map are 
\[\begin{aligned}
0=& h_2(y_1,y_2)\\
0=& h_{3k}(y_1,y_2)+x^\lambda
\end{aligned}\]
for $D^2(f)$. The space $D^3(f) \subset \CC^4$ is given by the following divided differences: 
\[\begin{aligned}
0=& y_1+y_2+y_3\\
0=& h_2(y_1,y_2,y_3)\\
0=& h_{3k-1}(y_1,y_2,y_3)\\
0=& h_{3k}(y_1,y_2)+x^\lambda.\\
\end{aligned}\]

We show that $D^3(f)$ does not have dimension zero. We do it by finding a curve in $D^3(f)$. Let $\zeta_3$ be a primitive third root of unity and choose $c \in \mathbb{C}$ such that $c^\lambda = -h_{3k}(1, \zeta_3)$. Define the parametrized analytic curve $\gamma : (\mathbb{C}, 0) \to (\mathbb{C}^4, 0)$ by
\[\gamma(t) = \left( x(t), y_1(t), y_2(t), y_3(t) \right) = \left( c \cdot t^{3k}, \; t^\lambda, \; t^\lambda \zeta_3, \; t^\lambda \zeta_3^2 \right).\]

Evaluating the defining equations of $D^3(f)$ along $\gamma(u)$, it follows by \cref{thm:sympolZEROS} that:
\[\begin{aligned}
& y_1(t) + y_2(t) + y_3(t) = t^{\lambda}(1 + \zeta_3 + \zeta_3^2) = 0\\
& h_2(y_1(t),y_2(t),y_3(t))=t^{2\lambda} h_2(1, \zeta_3, \zeta_3^2)=0\\
& h_{3k-1}(y_1(t), y_2(t), y_3(t)) = t^{\lambda(3k-1)} h_{3k-1}(1, \zeta_3, \zeta_3^2) = 0\\
& x(t)^\lambda + h_{3k}(y_1(t), y_2(t)) = (c t^{3k})^\lambda + t^{3k\lambda} h_{3k}(1, \zeta_3) = t^{3k\lambda}(c^\lambda + h_{3k}(1, \zeta_3)) = 0,\\
\end{aligned}\] where in third equation folows since $3k-1 \equiv 2 \pmod 3$.

Thus, $\gamma(t)$ is a non-constant analytic curve contained in $D^3(f)$ with $\gamma(0) = 0$. Therefore, $\dim D^3(f) \ge 1$, and by the Marar-Mond criterion \cref{thm: Marar-Mond crit}, $f$ is not $\mathcal{A}$-finite.
\end{proof}

One can iterate the previous process by adding several variables and increasing the dimension. More precisely, one can take a deformation of
\[g(y)=(y^a,y^b)\]
on several variables as follows:
\[f(x_1,\dots,x_{n-1},y)=\big(\xbm,y^a,y^b+m_1(\xbm,y)+\cdots+m_{n-1}(\xbm,y)\big),\] where $m_1,\ldots,m_{n-1} \in \mathbb{C}[x,y]$ are monomials.

We provide now the multiplicity two case.

\begin{proposition}\label{prop:mult2}
 Let $f : (\CC^{n},0)\to (\CC^{n+1},0)$ be a map germ given by
  \[f (\xbm, y)=\left(\xbm, y^2,y^{2k+1}+x_1^{\lambda_1}y+x_2^{\lambda_2}y+\ldots + x_{n-1}^{\lambda_{n-1}}y \right)\] for $k,\lambda_i \geq 1$ with $i\in \{1,\ldots,n-1\}$. Then, $f$ is $\Ascr$-finite and \[\codim_{\Ascr_e}(f)=k(\lambda_1-1)(\lambda_2-1)\cdots (\lambda_{n-1}-1).\]
  In particular, $f$ is $\Ascr$-stable if and only if $\lambda_j=1$ for some $j\in \{1,\ldots, n-1\}$.
 \end{proposition}
 \begin{proof}
 The first divided differences, that define $D^2(f)$, are
 \[\begin{aligned}
     0=& \ y_1+y_2,\\
     0=&\ \sum_{j+\ell=2k}y_1^jy_2^{\ell}+\sum_{i=1}^{n-1}x_i^{\lambda_i},
 \end{aligned}\]
 from which we see that there are no triple points. It is easy to see that $D^2(f)$ is isomorphic to
 \[0=y^{2k}+\sum_{\ell=1}^{n-1}x_i^{\lambda_i},\]
 which has isolated singularity, and it is smooth if some $\lambda_i$ is one. By the Marar-Mond criterion \cref{thm: Marar-Mond crit}, the germ $f$ is $\Ascr$-finitely determined.

 We compute now a monomial basis of $N\Ascr_e(f)$. Note that, in particular, we obtain in $f^*(\mathcal{O}_{n+1})$ the elements \[x_1^{s_1}x_2^{s_2}\cdots x_{n-1}^{s_{n-1}}y^{2\ell}, \ \ \mbox{for any} \ s_1,\dots,s_{n-1},\ell \geq 0.\] 
 Thus, we only need to analyze in $T\Ascr_e(f)$ the vectors whose components are monomials of the form $x_1^{s_1}\cdots x_{n-1}^{s_{n-1}}y^{2\ell+1}$, for any $s_1,\dots,s_{n-1},\ell \geq 0$.

 The Jacobian matrix of $f$ is given by 
 \[\left(\begin{array}{cccc}
1 & \cdots & 0 & 0\\
\vdots & \ddots & \vdots & \vdots\\
0 & \cdots & 1 & 0\\
0 & \cdots & 0 & 2y\\
\lambda_1x_1^{\lambda_1-1}y & \cdots & \lambda_{n-1}x_{n-1}^{\lambda_{n-1}-1}y  & (2k+1)y^{2k}+\sum_{i=1}^{n-1}x_i^{\lambda_i} \\
\end{array}\right).\]

By \cref{Tang-Space}, we obtain the following elements in $T\Ascr_e(f)$: 
\vspace{0.2cm}

\noindent 
(a) $\prod_i x_i^{s_i}\cdot y^{2\ell+1}e_j$ for each $j=1,\ldots, n-1$, and any $s_1,\dots,s_{n-1},\ell \geq 0$, by setting
\[
\begin{aligned}
\epsilon_j(\xbm,y) &= x_1^{s_1}\cdots x_{n-1}^{s_{n-1}}y^{2\ell+1},
& \epsilon_t &= 0 \ \text{for any } t\neq j,\\
\eta_{n+1}(f) &= -\lambda_jx_1^{s_1}\cdots x_j^{\lambda_j+s_j-1}\cdots x_{n-1}^{s_{n-1}}y^{2\ell+2},
& \eta_t &= 0 \ \text{for any } t\neq n+1.
\end{aligned}
\]

\noindent 
(b) $\prod_i x_i^{s_i}\cdot y^{2\ell+1}e_n$ for any $s_1,\dots,s_{n-1},\ell \geq 0$, by setting  
\begingroup
\newcommand{\Ms}{\prod_{i=1}^{n-1}x_i^{s_i}}
\[
\begin{aligned}
\epsilon_{n}(\xbm,y) &= \frac{1}{2}x_1^{s_1}\cdots x_{n-1}^{s_{n-1}} y^{2\ell},
& \epsilon_t &= 0 \ \text{for any } t\neq n,\\
\eta_{n+1}(f) &= -\frac{2k+1}{2}\Ms\,y^{2(k+\ell)}
-\frac{1}{2}\Ms\,y^{2\ell}\sum_{i=1}^{n-1}x_i^{\lambda_i},
& \eta_t &= 0 \ \text{for any } t\neq n+1.
\end{aligned}
\]
\endgroup

\noindent 
(c) $\prod_i x_i^{s_i}\cdot y^{2\ell+1}e_{n+1}$ for any $s_1,\dots,s_{n-1},\ell \geq 0$ and some $s_j \geq \lambda_j-1$, by setting
\[
\begin{aligned}
\epsilon_n(\xbm,y) &= \frac{1}{\lambda_j}x_1^{s_1}\cdots x_j^{s_j-\lambda_j+1}\cdots x_{n-1}^{s_{n-1}}y^{2\ell},
&\quad \epsilon_t &= 0 \ \text{for any } t\neq n,\\
\eta_{j}(f) &= -\frac{1}{\lambda_j}x_1^{s_1}\cdots x_j^{s_j-\lambda_j+1}\cdots x_{n-1}^{s_{n-1}}y^{2\ell},
&\quad \eta_t &= 0 \ \text{for any } t\neq j
\end{aligned}
\]

\noindent 
(d) $\ y^{2\ell+1}e_{n+1}$ for any $\ell \geq k$. In this case, by setting
\[
\begin{aligned}
\epsilon_n(\xbm,y) &= \frac{1}{2k+1}y^{2(\ell-k)+1},
&\quad \epsilon_t &= 0 \ \text{for any } t\neq n,\\
\eta_n(f) &= -\frac{2}{2k+1}y^{2(\ell-k+1)},
&\quad \eta_t &= 0 \ \text{for any } t\neq n
\end{aligned}
\]
we obtain that
\[
\left(y^{2\ell+1}+\frac{1}{2k+1}\sum_{i=1}^{n-1}x_i^{\lambda_i}y^{2(\ell-k)+1}\right)e_{n+1}\in T\Ascr_e(f).
\]
By item (c), the element
\[
\frac{1}{2k+1}\sum_{i=1}^{n-1}x_i^{\lambda_i}\cdot y^{2(\ell-k)+1}e_{n+1}\in T\Ascr_e(f),
\]
so we obtain $y^{2\ell+1}e_{n+1}$ as desired.
\newline

Then, a monomial basis of $N\Ascr_e(f)$ is given by the classes of \[B = \left\{ \begin{pmatrix} 0 \\ \vdots \\ 0 \\ \xbm^r y^{2\ell+1} \end{pmatrix} : \begin{array}{l} 0 \le \ell \le k-1, \\  0 \le r_j \le \lambda_j - 2 \text{ for each } j=1,\ldots, n-1 \end{array} \right\},\] where $\xbm^r = x_1^{r_1} x_2^{r_2} \cdots x_{n-1}^{r_{n-1}}$ for the exponent vector $r = (r_1, \dots, r_{n-1})$.

\vspace{0.2cm}
Thus, the $\Ascr_e$-codimension of $f$ is clearly \[\codim_{\Ascr_e}(f)=k(\lambda_1-1)(\lambda_2-1)\cdots (\lambda_{n-1}-1).\]
 \end{proof}

We show now that the multiplicity three differs from the previous case of multiplicity two.

\begin{proposition}
The map germ $f : (\CC^{n},0)\to (\CC^{n+1},0)$ given by
  \[f (\xbm, y)=\left(\xbm, y^3,y^{3k+2}+x_1^{\lambda_1}y+x_2^{\lambda_2}y+\ldots + x_{n-1}^{\lambda_{n-1}}y \right)\] is not $\Ascr$-finite if $n>2$.
\end{proposition}
\begin{proof}
Assume $n > 2$. Consider the complex subvariety $V \subset (\mathbb{C}^n, 0)$ defined by the equations
\[y = 0 \quad \text{and} \quad \sum_{\ell=1}^{n-1} x_{\ell}^{\lambda_{\ell}} = 0.\]
Since $n > 2$, the second condition defines a hypersurface in $\mathbb{C}^{n-1}$ of dimension $n-2 \ge 1$. Consequently, $V$ is a subvariety of positive dimension passing through the origin.

The Jacobian matrix of $f$ is:
 \[\begin{bmatrix}
 1 & \cdots & 0 & 0\\
 \vdots & \ddots & \vdots & \vdots\\
 0 & \cdots & 1 & 0\\
 0 & \cdots & 0 & 3y^2\\
 \lambda_1x_1^{\lambda_1-1}y & \cdots & \lambda_{n-1}x_{n-1}^{\lambda_{n-1}-1}y  & (3k+2)y^{3k+1}+\sum_{\ell=1}^{n-1}x_{\ell}^{\lambda_{\ell}} \\
 \end{bmatrix}\]
 Evaluating $Jf$ at any point $p = (\xbm, 0) \in V$, the last column vanishes completely. Thus, $\mathrm{rank}\, Jf(p) = n-1$, showing that $f$ has corank one at every point of $V$.

However, at any of those points away from the origin, it is not $\Ascr$-equivalent to any of the corank one stable singularities shown in \cref{thm:corankonestables}. Therefore, it has non-isolated instability and it is not $\Ascr$-finite by the Mather-Gaffney criterion of finite determinacy, see \cite[Theorem 4.5]{Mond2020}.
\end{proof}

\begin{remark}
    Observe that, considering $F:(\CC\times \CC,0)\to (\CC \times \CC^2,0)$ the one parameter stable unfolding
    \[F(t,y)=(t,y^2,y^{2k+1}+ty)\]
    of $(y^2,y^{2k+1})$, the map germ $f$ in \cref{prop:mult2} is the pullback of $F$ by $\gamma:(\CC^{n-1},0)\to (\CC,0)$ given by $\gamma(\xbm)=\textstyle\sum_{\ell=1}^{n-1}x_{\ell}^{\lambda_{\ell}}$, that is,
    \[\begin{array}{rcl}
\gamma^*F(x,y)&=&\big(\xbm,y^2,y^{2k+1}+\gamma(\xbm)y\big) \vspace{0.1cm}\\
&=&(\xbm,y^2,y^{2k+1}+x_1^{\lambda_1}y+x_2^{\lambda_2}y+\ldots + x_{n-1}^{\lambda_{n-1}}y ).\end{array}\]
    This is called an \textit{augmentation}, and it is well understood when augmentations are $\Ascr$-finite (see \cite[Theorem 3.3]{BrevaRibes2025}) and even, in some cases, what the generators of $N\Ascr_e f$ are, see \cite[Corollary 4.2]{BrevaRibes2022}. 
    Moreover, $f$ is a fold map, whose $\Ascr$-finiteness was characterized in terms of the $\Kscr$-finiteness of an auxiliary map in \cite[Proposition 2.3]{Houston1998} and whose $\Ascr_e$-codimension is determined using the same auxiliary map in \cite[Theorem 2.8]{Houston1998}.

    In the case of the multiplicity 3 germ $f$ in \cref{prop:mult3}, the one parameter unfolding
    \[F(t,y)=(t,y^3,y^{3k+2}+ty)\]
is not stable, so the same argument does not apply. Usually, augmentations of non-stable unfoldings are not $\Ascr$-finite.
\end{remark}

While for $n>2$ the map germ in the previous result fails to be $\mathcal{A}$-finite, the situation changes when $n=2$. In this case, the germ remains $\mathcal{A}$-finite, and we can determine its exact $\mathcal{A}_e$-codimension.
 \begin{proposition}
    \label{prop:mult3}
    The map germ $f : (\CC^{2},0)\to (\CC^{3},0)$ given by
  \[f (x, y)=\left(x, y^3,y^{3k+2}+x^{\lambda}y\right)\] is $\Ascr$-finite and satisfies
  \[\codim_{\Ascr_e}(f)=k(4\lambda - 3)+2\lambda -1,\] for all $k,\lambda \geq 1$.
 \end{proposition}
  \begin{proof}
 The idea of this proof is to provide a monomial basis of $N\Ascr_e(f)$. In $f^*(\Ocal_3)$ we have the elements
 $$x^ry^{3s}(y^{3k+2}+x^{\lambda}y)^t, \ \mbox{for any} \ r,s,t \geq 0.$$
Thus, we only need to consider the elements in $T\Ascr_e(f)$ whose components are monomials of the form $x^ry^{3s+1}$ and $x^ry^{3s+2}$ for any $r,s \geq 0$.
 
 By \cref{Tang-Space}, the extended tangent space to the $\Ascr$-orbit of $f$ is given by the following.
\[
   T\Ascr_e(f)=\left\{\left(\begin{array}{cc}
1 & 0\\
0 & 3y^2 \\
\lambda x^{\lambda -1}y & (3k+2)y^{3k+1}+x^{\lambda}\\
\end{array}\right) \left(\begin{array}{c}
\epsilon_1(x,y)\\
\epsilon_2(x,y)\\
\end{array}\right)+\left(\begin{array}{c}
\eta_1(f)\\
\eta_2(f)\\
\eta_3(f)
\end{array}\right); \begin{array}{c}
\epsilon \in \theta_2, \ \eta \in \theta_3,\\
\epsilon=(\epsilon_1,\epsilon_2), \\
\eta=(\eta_1,\eta_2,\eta_3)
\end{array}\right\}
\]

We denote the coordinates of $\CC^3$ by $X,Y$ and $Z$. The following vectors $v_1, \ldots, v_8$ belong to $T\Ascr_e(f)$:

\begin{align*}
v_1&=\big(x^ry^{3s+1},\,0,\,\lambda x^{r+\lambda-1}y^{3s+2}\big), && r,s\ge0,\\
&\qquad \epsilon=(x^ry^{3s+1},0),\ \eta=0;\\[3pt]
v_2&=\big(0,\,0,\,(3k+2)x^ry^{3s+2}+x^{r+\lambda}y^{3(s-k)+1}\big), && r\ge0,\ s\ge k,\\
&\qquad \epsilon=(0,x^ry^{3(s-k)+1}),\ \eta=(0,-3X^rY^{s-k+1},0);\\[3pt]
v_3&=\big(0,\,0,\,\tfrac{\lambda}{3k+2}\,x^{r+2\lambda-1}y^{3(s-k)+1}\big), && r\ge0,\ s\ge k,\\
&\qquad \epsilon=\big(\tfrac{1}{3k+2}x^{r+\lambda}y^{3(s-k)},0\big),\ \eta=\big(-\tfrac{1}{3k+2}X^{r+\lambda}Y^{s-k},0,0\big);\\[3pt]
v_4&=\big(x^ry^{3s+2},\,0,\,0\big), && r,s\ge0,\\
&\qquad \epsilon=(x^ry^{3s+2},0),\ \eta=(0,0,-\lambda X^{r+\lambda-1}Y^{s+1});\\[3pt]
v_5&=\big(0,\,3x^ry^{3s+1},\,x^{r+\lambda}y^{3s-1}\big), && r\ge0,\ s\ge1,\\
&\qquad \epsilon=(0,x^ry^{3s-1}),\ \eta=(0,0,-(3k+2)X^rY^{s+k});\\[3pt]
v_6&=\big(x^{r+1}y^{3s-2},\,0,\,\lambda x^{r+\lambda}y^{3s-1}\big), && r\ge0,\ s\ge1,\\
&\qquad \epsilon=(x^{r+1}y^{3s-2},0),\ \eta=0;\\[3pt]
v_7&=\big(0,\,3x^ry^{3s+2},\,(3k+2)x^ry^{3(s+k)+1}\big), && r,s\ge0,\\
&\qquad \epsilon=(0,x^ry^{3s}),\ \eta=(0,0,-X^{r+\lambda}Y^s);\\[3pt]
v_8&=\big(0,\,0,\,\lambda x^ry^{3s+1}\big), && r\ge\lambda-1,\ s\ge0,\\
&\qquad \epsilon=(x^{r-\lambda+1}y^{3s},0),\ \eta=(-X^{r-\lambda+1}Y^s,0,0).
\end{align*}

Note that $(x^ry^{3s+2},0,0) \in T\Ascr_e(f)$ for any $r,s \geq 0$ by $v_4$. We denote by $v_2'$ the vector in $T\Ascr_e(f)$, replacing $r$ by $r+\lambda -1$ in $v_2$. Then
\begin{equation}\label{3s+1-e_1}
(x^ry^{3s+1},0,0) \in T\Ascr_e(f), \ \mbox{for any}\ r\geq 0 \ \mbox{and} \ s\geq k,
\end{equation} since $(x^ry^{3s+1},0,0)=v_1-\frac{\lambda}{3k+2}v_2'+v_3$. Moreover,
\begin{equation}\label{eq:11'}
(x^ry^{3s+1},0,0)\in T\Ascr_e f\quad\text{for any } r\geq\lambda \text{ and } s\geq 0.
\end{equation}
Indeed, taking $\epsilon=0$ and $\eta=(X^{r-\lambda}Y^sZ,0,0)$, we obtain $\big(x^{r-\lambda}y^{3(s+k)+2}+x^ry^{3s+1},0,0\big)\in T\Ascr_e f$, and the first summand belongs to $T\Ascr_e f$ by $v_4$.

Our purpose now is to prove that:

\noindent \quad 
(a) $(0,0,x^ry^{3s+1}) \in T\Ascr_e(f)$ for any $\left\{\begin{array}{l}
0 \leq r \leq \lambda -2 \ \mbox{and} \ s \geq 2k+1 \vspace{0.1cm}\\
 r\geq \lambda -1 \ \mbox{and} \ s\geq 0. \\
 \end{array}\right.$

\vspace{0.2cm}
(b) $(0,0,x^ry^{3s+2}) \in T\Ascr_e(f)$ for any $\left\{\begin{array}{l}
r\geq 0\ \mbox{and} \ s \geq k \vspace{0.1cm}\\
 r\geq 2\lambda -1 \ \mbox{and} \ s\geq 0. \\
 \end{array}\right.$

\vspace{0.2cm}
 In fact, by $v_8$, we obtain $(0,0,x^ry^{3s+1}) \in T\Ascr_e(f)$, for $r\geq \lambda -1$ and $s\geq 0$. Moreover, $(0,0,x^ry^{3s+2}) \in T\Ascr_e(f)$, for any $r\geq 2\lambda -1$ and $s\geq 0$ since by the relation (\ref{eq:11'}), we conclude that
 \[(0,0,\lambda x^{r+\lambda -1}y^{3s+2})=v_1-(x^ry^{3s+1},0,0) \in T\Ascr_e(f),\qquad r\geq\lambda,\ s\geq0.\]

For $0 \leq r \leq \lambda -2$ in item (a), using the substitution $s \mapsto s-k$ in $v_7$, we obtain
$$w_7=(0,3x^ry^{3(s-k)+2},(3k+2)x^ry^{3s+1})\in T\Ascr_e(f), \ \mbox{for any} \ r\geq 0 \ \mbox{and}\ s\geq k.$$
Note that for $r,s\geq 0$ we have
\begin{equation}
    \label{3k+1-L}
x^ry^{3(s+k)+2}+x^{r+\lambda}y^{3s+1}=x^r(y^3)^{s}(y^{3k+2}+x^{\lambda}y) \in f^*(\Ocal_3).
\end{equation}
In (\ref{3k+1-L}), replacing power $s$ by $s-2k$, we obtain $x^ry^{3(s-k)+2}+x^{r+\lambda}y^{3(s-2k)+1} \in f^*(\Ocal_3)$ for $s \geq 2k$. Thus, taking $\epsilon=0$ and $\eta=(0,X^rY^{s-2k}Z,0)$, we obtain $u_7=(0,x^ry^{3(s-k)+2}+x^{r+\lambda}y^{3(s-2k)+1},0)\in T\Ascr_e(f)$, and consequently \[w_7':=(0, -3x^{r+\lambda}y^{3(s-2k)+1},(3k+2)x^ry^{3s+1})=w_7-3u_7\in T\Ascr_e(f),\] for any $r\geq 0$ and $s\geq 2k.$ Using the substitutions $r \mapsto r+\lambda$ and $s \mapsto s-2k$ in $v_5$, we obtain $w_5=(0,3x^{r+\lambda}y^{3(s-2k)+1}, x^{r+2\lambda}y^{3(s-2k)-1})\in T\Ascr_e(f)$. In this case, since $s\geq 1$ in $v_5$ we have $s\geq 2k+1$ in $w_5$. Then
\[w_8:=\big(0,\,0,\,(3k+2)x^ry^{3s+1}+x^{r+2\lambda}y^{3(s-2k)-1}\big)=w_7'+w_5 \in T\Ascr_e(f),\]
for any $r\geq 0$ and $s\geq 2k+1$. Since $x^{r+2\lambda}y^{3(s-2k)-1}=x^{r+2\lambda}y^{3(s-2k-1)+2}$ and $r+2\lambda\geq2\lambda-1$, the vector $(0,0,x^{r+2\lambda}y^{3(s-2k)-1})$ belongs to $T\Ascr_e f$ by (b).

Consequently, for any $1 \leq r\leq \lambda -2$ the element
\[(0,0,(3k+2)x^ry^{3s+1})\in T\Ascr_e(f), \ \text{for all}\ s\geq 2k+1.\]

It remains to analyze the case $(0,0,y^{3s+1})$ which will be addressed later. 

In (\ref{3k+1-L}), replacing $s$ by $s-k$ we have \begin{equation}
\label{Elem-f-O3}
x^ry^{3s+2}+x^{r+\lambda}y^{3(s-k)+1}\in f^*(\Ocal_3), \ \text{ for any } r\geq 0 \text{ and } s \geq k.
\end{equation} 

In this step, we determine which vectors of the form $(0,0,x^ry^{3s+2})$ belong to $T\Ascr_e(f)$ for $0 \leq r \leq 2\lambda -2.$

 Note that, by (\ref{Elem-f-O3}) the element 
$$(0,0,(3k+1)x^ry^{3s+2})=v_2-(0,0,x^ry^{3s+2}+x^{r+\lambda}y^{3(s-k)+1}) \in T\Ascr_e(f),$$ for any $0 \leq r\leq 2\lambda -2$ and $s\geq k$. Together with the case $r\geq2\lambda-1$ proved above, this completes the proof of (b).

By (\ref{3s+1-e_1}) and (\ref{eq:11'}), we obtain the elements $(x^ry^{3s+1},0,0)\in T\Ascr_e f  $ for $r\geq0$ and $s\geq k$, and for $r\geq\lambda$ and $s\geq0$. For $0\leq r\leq\lambda-1$ and $0\leq s\leq k-1$, we use $v_1$ to select the vector $(0,0,x^{r+\lambda-1}y^{3s+2})$ in the basis $B$ of the normal space $ N\Ascr_e f  $. Note that $r+\lambda-1\leq2\lambda-2$.

Lastly, we have $(0,x^ry^{3s+1},0)\in T\Ascr_e f  $ for $r\geq\lambda-1$ and $s\geq1$, by $v_5$ and (b): the third component of $v_5$ is $x^{r+\lambda}y^{3(s-1)+2}$, and $r+\lambda\geq2\lambda-1$. For the same reason, $(0,x^ry^{3s+1},0)\in T\Ascr_e f  $ for $0\leq r\leq\lambda-2$ and $s\geq k+1$. For $0\leq r\leq\lambda-2$ and $1\leq s\leq k$, $v_5$ gives
\[
(0,x^ry^{3s+1},0)\equiv-\tfrac13\,\big(0,0,x^{r+\lambda}y^{3(s-1)+2}\big)\pmod{ T\Ascr_e f  },
\]
which is an element of $B$, since $r+\lambda\leq2\lambda-2$ and $s-1\leq k-1$. Finally, for $s=0$ and $r\geq\lambda$, taking $\epsilon=0$ and $\eta=(0,X^{r-\lambda}Z,0)$, we obtain $(0,x^{r-\lambda}y^{3k+2}+x^ry,0)\in T\Ascr_e f  $, so
\[
(0,x^ry,0)\equiv-\big(0,x^{r-\lambda}y^{3k+2},0\big)\pmod{ T\Ascr_e f  },
\]
and the right-hand side is treated in the last paragraph below. The vectors $(0,x^ry,0)$ with $0\leq r\leq\lambda-1$ are elements of $B$.


Now, we turn to the case $(0,0,y^{3s+1})$. By $v_7$ and (\ref{Elem-f-O3}) for $r=0$ we obtain
 \[ \begin{array}{rcl}(0,-3x^{\lambda}y^{3(s-k)+1},(3k+2)y^{3(s+k)+1})&=&(0,3y^{3s+2},(3k+2)y^{3(s+k)+1}) \vspace{0.2cm}\\
 && -3(0,y^{3s+2}+x^{\lambda}y^{3(s-k)+1},0) \end{array}\] in $T\Ascr_e(f)$, for $s\geq k$. Then, $(0,0,y^{3s+1}) \in T\Ascr_e(f)$ for $s \geq 2k+1$, since by the previous argument $(0,x^ry^{3s+1},0)\in T\Ascr_e(f)$ for any $r\geq \lambda$ and $s\geq 1$, which we apply with $s-k\geq1$.
 
 Finally, by $v_7$ and (a) we have that $(0,x^ry^{3s+2},0)\in T\Ascr_e(f)$ for $r\geq \lambda-1$ and $s\geq 0$ and for $0\leq r\leq\lambda-2$ and $s\geq k+1$
 For any $0\leq r\leq\lambda-2$ and $0\leq s\leq k$, we use $v_7$ to select $(0,0,x^ry^{3(s+k)+1})$ in the basis $B$ of the normal space $ N\Ascr_e f  $.

Hence, a similar but shorter computations show independence of the following monomial basis of $N\Ascr_e(f)$,  given by the classes of \[B = \left\{   \begin{pmatrix} 0 \\ x^{r_1}y \\ 0 \\ \end{pmatrix}, \begin{pmatrix} 0 \\ 0 \\ x^{r_2}y^{3s+1} \\ \end{pmatrix}, \begin{pmatrix} 0 \\ 0 \\ x^{r_3}y^{3l+2}\\ \end{pmatrix} : \begin{array}{c}  0 \le r_1 \le \lambda -1, \\ \ 0 \le r_2 \le \lambda -2, \ \ 0 \le s \le 2k,  \\ 0 \le r_3 \le 2\lambda -2, \ \ 0 \le \ell \le k-1\end{array} \right\}.\]
Therefore, \[\codim_{\Ascr_e}(f)=\lambda +(\lambda -1)(2k+1)+(2\lambda -1)k=k(4\lambda-3)+2\lambda -1.\qedhere\]
 \end{proof}

 \vspace{1cm}

 \textbf{AI statement:} Generative AI was used exclusively for the layout of the equations, such as in the proof of Proposition 7.7. All the mathematical content was not obtained in any way or form from AI.

\bibliographystyle{myalpha.bst}
\bibliography{FullBib}
\end{document}